\documentclass[francais]{bourbaki}
 \usepackage[frenchb]{babel}
 \usepackage[T1]{fontenc}

 \usepackage{xspace,amssymb,mathrsfs,url,lmodern}

\title[Compter (rapidement)...]{Compter (rapidement) le nombre de solutions \\ 
       d'\'equations dans les corps finis}
\author{Antoine \MakeUppercase{Chambert-Loir}}
\date{Novembre 2006}
\bbkannee{59\textsuperscript{e} ann\'ee, 2006--2007}
\bbknumero{968}

\begin{abstract}
Le nombre de solutions dans les corps finis
d'un syst\`eme d'\'equations polynomiales
ob\'eit  \`a une tr\`es forte r\'egularit\'e, 
refl\'et\'ee par exemple par la rationalit\'e
de la fonction z\^eta d'une vari\'et\'e alg\'ebrique sur un corps fini,
ou la modularit\'e de la fonction $L$ de Hasse-Weil d'une courbe
elliptique sur $\Q$.
   Depuis une vingtaine d'ann\'ees des m\'ethodes efficaces
 ont \'et\'e invent\'ees pour calculer effectivement
ce nombre de solutions, notamment en vue d'applications
 \`a la cryptographie.
   L'expos\'e en pr\'esentera quelques-unes, g\'en\'eralement fond\'ees
l'utilisation de la formule des traces de Lefschetz dans
une th\'eorie cohomologique convenable,
  et expliquera leurs avantages respectifs. 
\end{abstract}

\begin{altabstract}
The number of solutions in finite fields of a system of polynomial
equations obeys a very strong regularity, reflected for example
by the rationality of the zeta function of an algebraic variety
defined over a finite field, or the modularity of Hasse-Weil's
$L$-function of an elliptic curve over~$\Q$.
Since two decades, efficient methods have been invented
to compute effectively this number of solutions, notably
in view of cryptographic applications.
This \emph{expos\'e} presents some of these methods,
generally relying on the use of Lefshetz's trace formula
in an adequate cohomology theory and discusses their
respective advantages.
\end{altabstract}

\altkeywords{algorithmic number theory, zeta functions, elliptic curves,
modular forms, $p$-adic cohomology} 

\alttitle{Counting (quickly) the number of solutions of equations
in finite fields}

\address{Universit\'e de Rennes~1 \\ 
Irmar (UMR 6625 du CNRS)  \\ Campus de Beaulieu \\
35042 Rennes Cedex}
\email{antoine.chambert-loir@univ-rennes1.fr}

\def\theenumi{\alph{enumi}}\def\labelenumi{\theenumi)}
\makeatletter
  \let\c@equation\c@defi
  \let\theequation\thedefi
  \let\cl@equation\cl@defi
  \numberwithin{paragraph}{subsection}
\makeatother

\let\bar\overline\let\tilde\widetilde
\def\Card#1{\mathopen| #1\mathclose|}
\let\card\Card
\def\dx{\mathrm dx}
\def\dt{\mathrm dt}
\def\A{\mathbf{A}}
\def\P{\mathbf{P}}
\def\F{\mathbf{F}}
\def\Z{\mathbf{Z}}
\def\C{\mathbf{C}}
\def\R{\mathbf{R}}
\def\Q{\mathbf{Q}}
\def\pgcd{\operatorname{pgcd}}

\def\Tr{\operatorname{Tr}}

\def\GL{\operatorname{GL}}
\def\SL{\operatorname{SL}}

\let\eps\varepsilon
\let\phi\varphi
\let\leq\leqslant\let\geq\geqslant

\def\cf{\emph{cf.}~}

\let\ra\rightarrow
\let\hra\hookrightarrow

\def\gm{\mathrm{G_m}}
\def\MW{{\text{\upshape MW}}}
\def\dR{{\text{\upshape dR}}}
\def\rig{{\text{\upshape rig}}}

\def\Gal{\operatorname{Gal}}
\def\Pic{\operatorname{Pic}}

\def\End{\operatorname{End}}
\def\id{\operatorname{id}}

\def\abs#1{\mathopen|{#1}\mathclose|}
\def\ord{\operatorname{ord}}
\def\rig{{\text{\upshape rig}}}

\def\cris{{\text{\upshape cris}}}

\def\legendre{\genfrac (){}0}
\def\Otilde{\widetilde{\mathrm O}}

\begin{document}
\maketitle

\setcounter{tocdepth}2
\tableofcontents

\section*{Introduction}

Soit $\F$ un corps fini et soit $f_1,\dots,f_m$
des polyn\^omes \`a coefficients dans~$\F$ en~$n$ ind\'etermin\'ees
$x_1,\dots,x_n$. 
Le but de cet expos\'e est de d\'ecrire des algorithmes 
permettant de calculer efficacement le nombre de solutions dans~$\F^n$
du syst\`eme d'\'equations $f_1=\dots=f_m=0$.

Pour tout entier~$k\geq 1$, notons $N_k$ le nombre
de solutions de ce syst\`eme dont les coordonn\'ees
appartiennent  au corps~$\F^{(k)}$, unique extension de~$\F$
de degr\'e~$k$ contenue dans une cl\^oture alg\'ebrique fix\'ee de~$\F$;
si $X$ est le sous-sch\'ema de l'espace affine 
d\'efini par l'annulation des~$f_i$,
on a donc $N_k=\Card{X(\F^{(k)})}$.
La fonction z\^eta~$Z(X,t)$ du sch\'ema~$X$ 
est alors donn\'ee par la formule
\begin{equation}
Z(X,t)=\exp \big( \sum_{k=1}^\infty N_k \frac{t^k}k \big)
\end{equation} 
et, ainsi que l'a d\'emontr\'e \textsc{Dwork}~\cite{dwork1960},
c'est une fraction rationnelle. Par cons\'equent, la suite~$(N_k)$ 
est d\'etermin\'ee par un nombre fini de ses termes.
Nous verrons aussi que ces algorithmes permettent de calculer
la fonction z\^eta $Z(X,t)$.

Tous les algorithmes d\'ecrits ci-dessous reposent 
sur un premier principe: il suffit, pour calculer $\Card{X(\F)}$,
de calculer une congruence $\Card{X(\F)}\equiv c\pmod N$,
o\`u $N$ est un entier strictement sup\'erieur \`a~$\Card{X(\F)}$,
par exemple $N>\Card\F^n$.
Plus g\'en\'eralement, il suffit que l'on connaisse 
un encadrement de~$\Card{X(\F)}$
de largeur inf\'erieure \`a~$N$;
c'est l\`a qu'interviendra l'analogue de l'hypoth\`ese de Riemann 
sur les corps finis, que \textsc{Deligne}~\cite{deligne74}
a d\'emontr\'ee, g\'en\'eralisant ainsi des r\'esultats
de~\textsc{Hasse} (courbes elliptiques) et \textsc{Weil} (courbes,
vari\'et\'es ab\'eliennes,...). Fixons-nous un tel encadrement
$C\leq \Card{X(\F)}<C+R$.

L\`a o\`u ces algorithmes diff\`erent, c'est sur la fa\c{c}on
de choisir un tel entier~$N$ puis de calculer~$c$.

Les premiers algorithmes, que nous qualifierons de~$\ell$-adiques,
font l'objet du chapitre~\ref{sec.ell} de ce rapport.
Ils ont pour arch\'etype l'algorithme d\'ecouvert en 1985
par R.~\textsc{Schoof}~\cite{schoof1985}
pour calculer le nombre de points d'une courbe
elliptique sur un corps fini.
Ces algorithmes choisissent un ensemble fini $\{\ell_1,\dots,\ell_s\}$
de {\og petits\fg} nombres premiers dont le produit~$L=\ell_1\dots\ell_s$
v\'erifie $L>R$ et calculent, pour tout~$i$ un \'el\'ement
$c_i\in\{0,\dots,\ell_i-1\}$ tel que $\Card{X(\F)}\equiv c_i\pmod{\ell_i}$.
Le th\'eor\`eme chinois permet d'en d\'eduire un entier~$c\in\{0,\dots,L-1\}$
tel que $\Card{X(\F)}\equiv c\pmod L$.
L'origine de la terminologie {\og$\ell$-adique\fg}
vient de ce qu'on peut interpr\'eter la congruence modulo~$\ell_i$ par
le calcul de la cohomologie \'etale modulo~$\ell_i$.

Hors du degr\'e~$1$ ou de ce qui en provient,
la cohomologie \'etale semble peu accessible au calcul
formel ; m\^eme dans ce cas, son calcul effectif
am\`ene rapidement \`a la consid\'eration de polyn\^omes de
tr\`es grand degr\'e. 
Le champ d'application
des algorithmes $\ell$-adiques est ainsi limit\'e aux courbes
de petit genre, aux vari\'et\'es ab\'eliennes de petite dimension.

N\'eanmoins, ces algorithmes sont polynomiaux
en le logarithme du cardinal de~$\F$: aussi bien le temps
de calcul que l'espace requis par le calcul sont major\'es
par une puissance de~$\log\Card\F$.

Nous pr\'esenterons au chapitre~\ref{sec.p} les algorithmes $p$-adiques,
o\`u l'entier~$p$ d\'esigne la caract\'eristique du corps~$\F$. 
Ils proc\`edent en effet en choisissant pour~$N$
une puissance de~$p$ 
et en calculant (plus ou moins)
la cohomologie $p$-adique de~$X$ modulo~$N$.
Par cohomologie $p$-adique, j'entends ici la 
cohomologie de Monsky-Washnitzer et ses avatars
(rigide, cristalline), qui sont des analogues de la cohomologie
de De Rham. D\'efinie comme cohomologie d'un complexe
explicite, la cohomologie $p$-adique se pr\^ete naturellement
bien au calcul effectif
et l'on peut esp\'erer appliquer ces m\'ethodes dans des situations
g\'eom\'etriques tr\`es g\'en\'erales. Malgr\'e tout, il semble que seules les
courbes et les surfaces aient fait l'objet
d'impl\'ementations pouss\'ees.

Toutefois, parce qu'ils demandent de manipuler
des polyn\^omes de degr\'es au moins~$p$, la d\'ependance en~$\log p$
de leur complexit\'e n'est pas polynomiale.
Ils n'en restent pas moins des algorithmes de choix lorsque
$p$ est petit, notamment dans les applications cryptographiques
o\`u l'on a souvent $p=2$.

\bigskip

Au fur et \`a mesure du d\'eveloppement de ces algorithmes,
ils ont \'et\'e programm\'es et leurs performances \'eprouv\'ees
\`a l'aune des records qu'ils permirent d'obtenir,
c'est-\`a-dire le calcul de~$\Card{X(\F)}$
pour des corps~$\F$ de cardinal le plus grand possible.
Lorsque $X$ est une courbe elliptique, on a pu atteindre
un cardinal~$q$ de plus de~$2\,000$ chiffres (en base~$10$)
par l'algorithme de~\textsc{Schoof}, et d'environ~$40\,000$
chiffres (mais en caract\'eristique~$p=2$) 
par l'algorithme $2$-adique de~\textsc{Mestre}. Ces
calculs ont pris plusieurs mois.
La diminution de l'espace m\'emoire n\'ecessit\'e par ces algorithmes
a aussi fait l'objet de travaux importants.

Parall\`element, ils ont trouv\'e un champ d'application dans
la cryptographie \`a clef publique et se sont retrouv\'es
au c\oe ur de logiciels commerciaux. Comme nous le verrons
plus bas, les corps~$\F$ qu'il faut alors manipuler sont
de taille bien plus modeste, disons une cinquantaine de chiffres d\'ecimaux.

\medskip

Le premier chapitre de ce texte est consacr\'e \`a quelques applications
de ce probl\`eme algorithmique et de ses diverses solutions efficaces.
J'exposerai ensuite les grandes lignes
de la plupart des algorithmes $\ell$-adiques, puis $p$-adiques, 
actuellement utilis\'es. 
Il s'av\`ere en fait qu'une bonne partie de la th\'eorie g\'en\'erale 
et abstraite d\'evelopp\'ee au \textsc{xx}\textsuperscript e
si\`ecle dans l'\'etude des conjectures de Weil
donne naturellement lieu \`a des algorithmes efficaces.
Cependant, cette constatation n'est pas all\'ee de soi 
et le cr\'edit en revient bien aux math\'ematiciens
tels que~\textsc{Schoof}, \textsc{Elkies}, \textsc{Atkin}
(pour la partie $\ell$-adique), \textsc{Satoh}, \textsc{Mestre},
\textsc{Kedlaya}, \textsc{Lauder} (pour la
partie~$p$-adique) dont les noms \'emailleront ce texte.
\`A moins d'achever cet expos\'e juste apr\`es
le chapitre consacr\'e aux applications,
il m'a ainsi fallu d\'epasser 
le lapidaire et spontan\'e {\og On peut le faire!\fg}
sans pour autant plonger le lecteur dans la complexit\'e
ph\'enom\'enale des id\'ees suppl\'ementaires qui ont \'et\'e n\'ecessaires
\`a l'obtention des records \'evoqu\'es plus haut.
Le compromis que j'ai essay\'e d'adopter dans ce texte,
un peu diff\'erent des nombreux survols du sujet disponibles dans
la litt\'erature, est celui d'un math\'ematicien pur subitement int\'eress\'e par
ce probl\`eme de math\'ematiques appliqu\'ees. 

\medskip

Lorsque je d\'ecris la complexit\'e d'algorithmes en temps
ou en espace, j'emploie
les notations~$\mathrm O(\cdot)$ et~$\Otilde(\cdot)$.
La premi\`ere signifie que le nombre d'op\'erations \'el\'ementaires,
resp. l'espace disque, requis par l'algorithme est
major\'e par un multiple de son argument, lorsque celui-ci
tend vers l'infini. La seconde est analogue, \`a une puissance
du logarithme de l'argument pr\`es;
en pratique, il suffit de retenir que
$\Otilde(x)$ est major\'e par $\mathrm O(x^{1+\eps})$
pour tout~$\eps>0$.
Toutefois, m\^eme si je n'en parlerai jamais,
il ne faut pas perdre de vue que le contr\^ole
de la constante que cachent ces notations 
est d'une importance pratique
capitale; il est bien diff\'erent de pouvoir
obtenir un r\'esultat en une minute plut\^ot qu'en mille.

\medskip

Je tiens \`a remercier Jean-Beno\^{\i}t~\textsc{Bost},
Bas~\textsc{Edixhoven}, 
Reynald \textsc{Lercier}, Bernard \textsc{Le Stum},
David \textsc{Lubicz}
et Jean-Fran\c{c}ois \textsc{Mestre}
de l'aide qu'ils m'ont apport\'ee
au cours de la pr\'eparation de cet expos\'e.
Je remercie aussi Robert \textsc{Carls}, David \textsc{Kohel}, 
Ren\'e \textsc{Schoof} et Jean-Pierre~\textsc{Serre} 
pour leurs commentaires sur la premi\`ere version de ce texte.

\section{Applications}

\subsection{Crit\`eres de primalit\'e}

\^Etre en mesure de d\'ecider si un entier naturel est ou pas un nombre
premier est une question arithm\'etique fondamentale
dont les techniques modernes de cryptographie ont d'ailleurs
accru l'importance.

En 1986, S.~\textsc{Goldwasser} et J.~\textsc{Kilian}
ont propos\'e (voir~\cite{goldwasser-kilian1999})
le premier algorithme permettant de d\'ecider
si un entier~$N$ est un nombre premier dont la complexit\'e 
soit polynomiale en~$\log N$. Cet algorithme requiert
de calculer le cardinal de courbes elliptiques~$E$ sur
l'anneau $\Z/N\Z$ {\og choisies au hasard\fg}. Pour cela,
on peut tenter d'appliquer l'algorithme de \textsc{Schoof}, en faisant
comme si $N$ \'etait premier. Si l'algorithme \'echoue, cela prouve
que~$N$ n'est pas premier. Supposons qu'il fournisse un cardinal putatif~$c$.
On peut tester si un point au hasard~$P$ sur la courbe~$E$ est
annul\'e  par~$c$; si ce n'est pas le cas, $N$ n'est pas premier.
Inversement, supposons que l'ordre~$d$ de~$P$  poss\`ede
un facteur premier~$p$ tel que $p>(1+\sqrt N)^2$
et tel que le point~$[d/p]P$ ne rencontre pas l'origine~$O$ de la
courbe~$E$ (au sens o\`u ce point~$[d/P]P$ ait des coordonn\'ees
homog\`enes~$(x:y:z)$ dans~$(\Z/N\Z)$, $z$ \'etant premier \`a~$N$);
alors $N$ est premier. (Sinon, d\'esignant par~$\ell$ le plus petit
facteur premier de~$N$, l'image de~$P$
dans~$E(\Z/\ell\Z)$ serait d'ordre multiple de~$p$,
et cela contredirait la borne de Hasse pour le cardinal 
d'une courbe elliptique sur un corps fini). L'algorithme
de~\textsc{Goldwasser} et~\textsc{Kilian} tente
alors d'exhiber de telles familles~$(E,P,d,p)$
o\`u $\Card{E(\Z/N\Z)}=2p$,
la primalit\'e de~$p$ \'etant \'etablie r\'ecursivement par la m\^eme m\'ethode.

Comme l'algorithme de~\textsc{Schoof} est de complexit\'e
polynomiale en~$\log N$, il en est de m\^eme de
celle de l'algorithme de~\textsc{Goldwasser} et~\textsc{Kilian}.
Toutefois, le fait que cet algorithme parvienne \`a conclure
pour tout~$N$ d\'epend d'une conjecture apparemment hors de port\'ee
sur la r\'epartition des nombres
premiers dans de petits intervalles.

\textsc{Adleman} et~\textsc{Huang}~\cite{adleman-huang1992} ont eu l'id\'ee
d'utiliser des courbes de genre~$2$.
Cela fournit plus de latitude et leur permet d'affirmer
l'existence d'un algorithme probabiliste
de complexit\'e polynomiale en~$\log N$
permettant de d\'ecider si l'entier~$N$  est premier.

Cependant, l'algorithme de~\textsc{Schoof}, m\^eme avec
les am\'eliorations d'\textsc{Atkin} et \textsc{Elkies}
(qui font l'objet du paragraphe~\ref{subsec.sea} ci-dessous)
n'est pas suffisamment efficace pour permettre d'envisager
de tester ainsi la primalit\'e d'entiers ayant plus de quelques
centaines de chiffres d\'ecimaux. 
Si l'algorithme~fast\textsc{ecpp}
(\emph{fast elliptic curve primality proving})
a permis de prouver la primalit\'e de nombres ayant
plus de vingt mille chiffres d\'ecimaux (\textsc{Morain}, mi-2006),
c'est en utilisant l'id\'ee
d'\textsc{Atkin} d'employer, plut\^ot que des courbes al\'eatoires,
des courbes elliptiques \`a multiplication complexe 
(de discriminants relativement petits, au plus $10^6$)
dont l'algorithme de \textsc{Cornacchia}~\cite{nitaj1995}
bien que probabiliste, permet de calculer rapidement le cardinal.

Pour plus de d\'etails,
je renvoie \`a l'expos\'e de F.~\textsc{Morain}~\cite{morain2004} \`a ce S\'eminaire,
au chapitre~9 du livre~\cite{cohen1993} d'H.~\textsc{Cohen},
ainsi qu'au chapitre~25 du manuel~\cite{cohen2006}.

\subsection{Calcul d'une racine carr\'ee modulo~$p$}

Soit $p$ un nombre premier impair et soit $a$ un \'el\'ement de~$\F_p$.
Le calcul du symbole de Legendre $\legendre ap=a^{(p-1)/2}$
est un moyen simple pour d\'ecider si $a$ est un carr\'e,
\`a condition bien s\^ur de calculer la puissance par des \'el\'evations au carr\'e
successives.
Ce que ce crit\`ere ne dit pas, c'est
comment trouver une racine carr\'ee, c'est-\`a-dire
un \'el\'ement $b\in\F_p$ tel que $b^2=a$.

Lorsque $p\equiv 3\pmod 4$, on peut poser $b=a^{(p+1)/4}$.
Le cas crucial est donc celui o\`u $p\equiv 1\pmod 4$,

Les algorithmes de factorisation de polyn\^omes dans~$\F_p$,
tel celui de Berlekamp, fournissent une solution
efficace.  Rappelons-en le principe dans ce cas particulier:
si $x$ est un \'el\'ement de~$\F_p$, distinct de~$\pm b$, 
les deux racines
$x+b$ et~$x-b$ du polyn\^ome $(X-x)^2-a$ seront simultan\'ement
carr\'es ou non carr\'es si et seulement si $(x+b)/(x-b)$ est un carr\'e
dans~$\F_p$. Si ce n'est pas le cas, c'est-\`a-dire une fois sur deux
si $x$ est choisi au hasard,
le pgcd des polyn\^omes~$X^2-a$ et $(X+x)^{(p-1)/2}-1$ sera 
l'un des polyn\^omes $X\pm b$. 
Le calcul d'un tel pgcd requiert un nombre d'op\'erations 
\'el\'ementaires au plus \'egal \`a une puissance de~$\log p$ : 
il suffit en effet de calculer $(X+x)^{(p-1)/2}-1$ 
dans la $\F_p$-alg\`ebre $\F_p[X]/(X^2-a)$.

Comme les autres algorithmes de factorisation dans les corps finis
de complexit\'e \'equivalente,
il s'agit toutefois d'un algorithme probabiliste
c'est-\`a-dire
que l'on a seulement une tr\`es forte probabilit\'e que l'algorithme
se termine en un temps donn\'e.
Pr\'ecis\'ement, la probabilit\'e que l'algorithme se termine en~$N$
\'etapes est $1-2^{-N}$, mais rien n'interdit de n'avoir pas de chance.

Connaissant un g\'en\'erateur du sous-groupe $2$-primaire
de~$\F_p^*$,
l'algorithme de~\textsc{Shanks} pr\'esent\'e dans~\cite{shanks1971}
permet alors de calculer des racines carr\'ees dans~$\F_p$
en temps~$\Otilde(\log p)$.
D'apr\`es~\cite{bach1990},
si l'hypoth\`ese de Riemann g\'en\'eralis\'ee aux fonctions~$L$
de Dirichlet est v\'erifi\'ee, le plus petit entier positif~$x$
qui n'est pas un carr\'e modulo~$p$
v\'erifie $x < 2(\log p)^{2}$.
\'Ecrivons $p-1$ sous la forme $2^kq$ avec $k\geq 0$ et~$q$ impair;
on voit que $\xi=x^{(p-1)/q} \mod p$ est un g\'en\'erateur 
du sous-groupe~$2$-primaire de~$\F_p^*$. 
Modulo~GRH, on a ainsi un algorithme pour calculer des racines carr\'ees
dans~$\F_p^*$ dont la complexit\'e est polynomiale en~$\log p$.
H\'elas, on ne sait pas en g\'en\'eral construire un tel g\'en\'erateur~$\xi$
de mani\`ere d\'eterministe en temps polynomial en~$\log p$
sans faire appel \`a l'hypoth\`ese de Riemann.

Comme l'a montr\'e~\textsc{Schoof},
une cons\'equence de son algorithme de calcul du cardinal
d'une courbe elliptique sur~$\F_p$ est un algorithme
d\'eterministe (d\'ependant de~$a$) de complexit\'e polynomiale
en~$\log p$ pour calculer une racine carr\'ee de~$a$ modulo~$p$.

Supposons en effet 
que $a \equiv -D\pmod p$, o\`u $-D$ est le discriminant d'un
ordre~$\mathscr O$ d'un corps quadratique imaginaire~$K$.
Par la th\'eorie de la multiplication complexe,
\textsc{Schoof} construit une courbe elliptique~$E$ sur
une extension finie~$\F$ de~$\F_p$ dont l'anneau des
endomorphismes est~$\mathscr O$.
C'est la partie de l'algorithme la plus co\^uteuse
car $[\F:\F_p]$ est de l'ordre de~$\sqrt {-D}$; la
complexit\'e de son algorithme d\'epend donc de l'entier~$a$.
Soit $q$ le cardinal de~$\F$.
Dans le corps~$K$, l'endomorphisme de Frobenius
est de la forme $x+y\sqrt D$, avec $x,y\in\frac12\Z$;
on a $2x=t_E=q+1-\Card{E(\F)}$ et $x^2+Dy^2=q$.
Une fois calcul\'e $\Card{E(\F)}$ par la m\'ethode de~\textsc{Schoof}
(et non par celle de~\textsc{Cornacchia}, probabiliste),
on conna\^{\i}t ainsi~$x$ et~$\abs y$.
Comme $D$ est un carr\'e modulo~$p$, la courbe elliptique $E$
est ordinaire et $t_E$ n'est pas multiple de~$p$; par
suite~$x$ et $y$  ne sont pas multiples de~$p$
et la r\'eduction modulo~$p$ de $(2x)/(2y)$ est une racine carr\'ee de~$a$.

Comme le remarque~\textsc{Schoof} \`a la fin de~\cite{schoof1985},
on peut combiner ceci
avec l'algorithme de~\textsc{Shanks} lorsqu'on suppose, par exemple,
$p\not\equiv 1\pmod{16}$. En effet, si $\xi$ est
un g\'en\'erateur de la composante
$2$-primaire de~$\F_p^*$, c'est une racine de l'unit\'e dont
l'ordre divise~$8$,  
 donc $\xi$ s'exprime en termes de racines carr\'ees de~$-1$ et~$2$
dans~$\F_p$. On peut ainsi calculer~$\xi$ en temps polynomial en~$\log p$
par l'algorithme de~\textsc{Schoof} utilisant les courbes elliptiques.

\medskip

De mani\`ere analogue, \textsc{Pila} a montr\'e dans~\cite{pila1990}
comment son algorithme de calcul de la fonction z\^eta
d'une courbe de genre sup\'erieur, appliqu\'e \`a la courbe
de Fermat d'\'equation $x^\ell+y^\ell+z^\ell=0$ 
o\`u $\ell$ est un nombre premier fix\'e, permet
de calculer en temps polynomial en~$\log p$
les racines primitives $\ell$-i\`emes de l'unit\'e dans
le corps fini~$\F_p$, pourvu bien s\^ur que $p\equiv 1\pmod \ell$.
\subsection{Cryptographie}

En 1975, \textsc{Diffie} et~\textsc{Hellman} 
ont propos\'e~\cite{diffie-hellman1976}
une solution
\'el\'egante permettant \`a deux individus d'\'echanger une information secr\`ete
bien que le canal de communication puisse \^etre espionn\'e par une
oreille indiscr\`ete.

Son principe est le suivant. Les deux protagonistes, Antoine
et Bernadette, conviennent d'un groupe~$G$ (not\'e multiplicativement)
et d'un \'el\'ement~$g$ de ce groupe.
Antoine choisit un entier~$a$, calcule~$g^a$ et le transmet \`a Bernadette;
celle-ci choisit un entier~$b$, calcule~$g^b$ et le transmet \`a Antoine.
Le secret commun est l'\'el\'ement $g^{ab}$ du groupe~$G$,
que nos deux h\'eros sont en mesure de calculer 
puisque $g^{ab}=(g^a)^b=(g^b)^a$;
ils peuvent par exemple l'utiliser 
comme param\`etre d'un syst\`eme de codage sym\'etrique.

Il est n\'ecessaire d'indiquer 
que ce protocole ne r\'esiste pas \`a une attaque active: 
supposons que Charles s'immisce dans la conversation 
et parvienne \`a se faire passer
pour Antoine \`a Bernadette et \`a Bernadette pour Antoine.
Il peut alors choisir des entiers~$a'$, $b'$, transmettre~$g^{a'}$
\`a Bernadette et $g^{b'}$ \`a Antoine. Ce dernier utilise
donc~$g^{ab'}$ pour coder ou d\'ecoder un message, tandis que Bernadette
utilise~$g^{a'b}$.  Puisqu'il conna\^{\i}t~$g^{ab'}=(g^a)^{b'}$
et $g^{a'b}=(g^b)^{a'}$, Charles peut intercepter un message d'Antoine,
le d\'ecoder et le recoder \`a l'intention de Bernadette, ou inversement,
sans qu'aucun des deux n'ait pu se douter de quoi que ce soit.  

M\^eme s'il n'a pu intervenir physiquement dans la conversation,
Charles a connaissance de~$G$, $g$
ainsi que des deux \'el\'ements $g^a$ et~$g^b$. Pour qu'il puisse
en d\'eduire~$g^{ab}$, il suffirait qu'il soit en mesure de
calculer~$a$ (ou~$b$).

Le probl\`eme, \'etant donn\'e deux \'el\'ements $g$ et~$h$
d'un groupe~$G$, de d\'eterminer un entier~$a$ tel que $h=g^a$
est appel\'e \emph{probl\`eme du logarithme discret}.
Pour que le protocole de Diffie--Hellman r\'esiste \`a une
attaque passive, il est manifestement
n\'ecessaire que le probl\`eme du logarithme discret dans le groupe~$G$
soit difficile \`a r\'esoudre en pratique;
voir~\cite{maurer-wolf1999,muzereau-smart-vercauteren2004}
pour l'\'etude de la r\'eciproque, conjecturalement vraie
--- il suffirait de savoir construire, pour tout facteur
premier~$p$ de~$\Card G$, une courbe elliptique sur~$\F_p$,
dont le nombre de points est {\og lisse\fg}, c'est-\`a-dire
que ses facteurs premiers sont petits.

Les groupes cycliques $\Z/n\Z$ ne conviennent \'evidemment pas,
car l'algorithme d'Euclide permet tr\`es facilement,
$g$ \'etant un g\'en\'erateur fix\'e de ce groupe, de calculer $a\pmod n$
connaissant $ga\pmod n$. 
\textsc{Diffie} et~\textsc{Hellman} ont propos\'e 
d'utiliser les groupes multiplicatifs de corps finis.

En 1985, \textsc{Koblitz} et~\textsc{Miller} ont montr\'e
que les groupes form\'es des points d'une courbe elliptique
sur un corps fini sont de bons candidats; plus g\'en\'eralement,
on peut imaginer utiliser les groupes~$\Pic^0 C$
des diviseurs de degr\'e~$0$ sur une courbe~$C$ d\'efinie 
sur un corps fini~$\F$.
Il s'agit toutefois de trouver un bon compromis entre
la commodit\'e du calcul 
et la difficult\'e du probl\`eme du logarithme discret
dans le groupe~$G$.

Si $G$ est d'ordre~$n$, il y a de nombreux algorithmes
en~$\Otilde(\sqrt n)$ pour r\'esoudre ce probl\`eme du logarithme discret,
citons celui des kangourous~\textsc{Pollard}~\cite{pollard1978}
qui repose sur le \emph{paradoxe des 
anniversaires} : tirons au hasard, avec remise, des \'el\'ements d'un ensemble
fini de cardinal~$n$ ; le nombre moyen de tirages avant
qu'on  obtienne un \'el\'ement d\'ej\`a tir\'e est 
$\sqrt{\pi n/2}$. 

Lorsque l'on conna\^{\i}t une factorisation de~$n$,
on peut tenter de r\'esoudre le logarithme discret dans chacun
des quotients d'ordre premier de~$G$, puis utiliser le th\'eor\`eme chinois,
d'o\`u un algorithme de complexit\'e~$\Otilde(\sqrt t)$
si $t$ est le plus grand facteur premier de~$n$
(algorithme de \textsc{Pohlig}-\textsc{Hellman}, \cite{pohlig-hellman1978}).

Dans un groupe {\og g\'en\'erique\fg}, c'est-\`a-dire dont
on ne sait rien et dont un oracle calcule
le produit de deux \'el\'ements, l'inverse d'un \'el\'ement et teste l'\'egalit\'e
de deux \'el\'ements, V.~\textsc{Shoup}~\cite{shoup1997}  a montr\'e 
qu'un algorithme requiert un nombre de recours
\`a l'oracle au moins proportionnel \`a~$\sqrt n$.

Compte tenu de ces attaques et de la puissance des moyens de calcul
actuels,
l'entier~$n$ doit donc \^etre au moins \'egal \`a~$2^{160}$,
de m\^eme que son plus grand facteur premier.
D'o\`u en particulier la n\'ecessit\'e de conna\^{\i}tre l'ordre du groupe~$G$,
donc de savoir calculer le cardinal
d'une courbe elliptique
ou, plus g\'en\'eralement, de la jacobienne d'une courbe
d\'efinie sur un corps fini.

Pour les groupes~$G=\F_q^*$, $G=\Pic^0(C)$,
mentionn\'es plus haut, il y a de nombreuses tentatives
pour calculer le logarithme discret.
Les techniques d'indice fournissent par exemple
des algorithmes sous-exponentiels
pour les groupes multiplicatifs de corps finis;
pour les groupes de classes de diviseurs de courbes hyperelliptiques,
ils sont plus efficaces que les algorithmes g\'en\'eriques
lorsque le genre est~$> 4$.
La non d\'eg\'enerescence de l'accouplement de Tate--Lichtenbaum
\[ \Pic^0(C)[n] \times\Pic^0(C_{\F_{q^k}})/n \ra \F_{q^k}^*/(\F_{q^k}^*)^n \]
permet une r\'eduction du probl\`eme au cas du groupe
multiplicatif d'une extension finie de~$\F_q$:
si $q^k\equiv 1\pmod n$, cet accouplement (compos\'e avec l'\'el\'evation
\`a la puissance $(q^k-1)/n$) fournit une injection de~$G$
dans~$\F_{q^k}^*$.
Comme cet accouplement se calcule ais\'ement~\cite{frey-ruck1994},
il convient donc de choisir des courbes~$C$ telles que 
l'ordre multiplicatif de~$q$ modulo~$n$ ne soit pas trop petit,
en pratique $\leq 20$. En particulier, cela proscrit
les courbes elliptiques supersinguli\`eres.
De m\^eme, il faut \'eviter que l'ordre~$n$ de~$\Pic^0(C)$ ne 
soit multiple  de la caract\'eristique du corps~$\F$.

Je renvoie aux ouvrages~\cite{blake-seroussi-smart1999,blake-seroussi-smart2005,cohen2006} pour une description d\'etaill\'ee des diverses attaques possibles, 
des choix raisonnables d'un corps fini, d'une courbe~$C$ et d'un \'el\'ement
de~$\Pic^0(C)$, et de la fa\c{c}on dont tout ceci peut
\^etre implant\'e dans une carte \`a puce.
En outre, les articles~\cite{koblitz-menezes2004,galbraith-menezes2005} m'ont 
\'et\'e tr\`es utiles pour \'ecrire ce paragraphe.

Signalons enfin que les protocoles cryptographiques
reposant sur la difficult\'e de r\'esoudre le logarithme
discret dans une courbe elliptique
ont fait l'objet d'une sp\'ecification par
divers organismes de normalisation (ANSI, ISO, etc.)
et sont au c\oe ur  de nombreux syst\`emes cryptographiques
commercialis\'es.

\section{L'approche $\ell$-adique}\label{sec.ell}

Soit $\F$ un corps fini \`a~$q$ \'el\'ements et soit $E$ une courbe
elliptique sur~$\F$, donn\'ee par une \'equation plane (inhomog\`ene)
de la forme 
\begin{equation}\label{eq.weierstrass.general}
  y^2+a_1xy+a_3y =x^3+a_2x^2+a_4 x+a_6,
\end{equation}
o\`u les coefficients $a_1,a_2,a_3,a_4,a_6$ sont des \'el\'ements de~$\F$.
Lorsque $p>3$, on peut se ramener \`a une \'equation de la forme
de Weierstrass:
\begin{equation}\label{eq.weierstrass}
y^2=x^3+ a x + b,
\end{equation}
o\`u $a$ et~$b$ sont des \'el\'ements de~$\F$.
Pour tout corps~$k$ contenant~$\F$, notons $E(k)$ les $k$-points de~$E$,
c'est-\`a-dire l'ensemble des solutions de l'\'equation~\eqref{eq.weierstrass.general}
auquel l'on adjoint le point \`a l'infini~$O$ de~$\P^2(k)$
de coordonn\'ees homog\`enes~$(0:1:0)$.
On sait que $E$ poss\`ede une unique structure de groupe
alg\'ebrique, commutatif, dont l'\'el\'ement neutre est le point~$O$.
Elle se d\'eduit de la construction par s\'ecantes et tangentes:
pour tout corps~$k$ contenant~$\F$,
trois points $P_1$, $P_2$ et $P_3$ de~$E(k)$ sont align\'es
si et seulement l'on a $P_1+P_2+P_3=O$ dans le groupe $E(k)$.

Nous voulons calculer le cardinal du groupe ab\'elien $E(\F)$.

\subsection{Premi\`eres approches}
\paragraph{Symboles de Legendre}
Lorsque $E$ est de la forme~\ref{eq.weierstrass},
on a 
\[  \Card {E(\F)} = 1+ \sum_{x\in\F} \begin{cases}
        2 & \text{si $x^3+ax+b\in(\F^*)^2$} \\
       1 & \text{si $x^3+ax+b=0$} \\
     0 &\text{sinon} \end{cases}
 = q+1+ \sum_{x\in\F}\legendre{x^3+ax+b}{q}, 
\]
o\`u $\legendre{t}p$ est le symbole de Legendre dans~$\F$,
qui vaut $1$ si $t$ est un carr\'e non nul dans~$\F$, $0$
si $t=0$ et $-1$ sinon.
La complexit\'e de cette m\'ethode est $\Otilde(q)$.
Elle n'est ainsi utilisable que si $q$ est petit:
un programme \'el\'ementaire 
requiert d\'ej\`a 10~s pour d\'eterminer que lorsque $q=1\,048\,609$,
et $E$ a pour \'equation $y^2=x^3-1$, $\Card {E(\F)}=1\,049\,412$.
 

\paragraph{Frobenius}
L'\'el\'evation des coordonn\'ees \`a la puissance~$q$
d\'efinit un morphisme $\pi_E\colon E\ra E$ de groupes alg\'ebriques,
appel\'e endomorphisme de Frobenius.
Dans l'anneau $\End(E)$ des endomorphismes de~$E$, $\pi_E$
v\'erifie une relation polynomiale:
\begin{equation}\label{eq.frobenius}
\pi_E^2 - t_E \pi_E + q = 0, 
\end{equation}
o\`u $t_E$ est un entier relatif tel que $\abs{t_E}\leq 2\sqrt q$.
En outre (Hasse, \cite{hasse1933,hasse1936}), 
cet entier $t_E$ est reli\'e au cardinal de~$E(\F)$
par la formule
\begin{equation}\label{eq.hasse}
\Card{E(\F)} = q+1-t_E.
\end{equation}
Rappelons enfin que $t_E$ est la trace de~$\pi_E$
dans~$\End(E)$;
on a en effet l'\'egalit\'e :
\begin{equation}
t_E \id_E = \pi_E+\pi_E^\vee, 
\end{equation} 
o\`u $\pi_E^\vee$ est l'isog\'enie duale~$\pi_E^\vee$,
d\'efinie par $\pi_E^\vee\circ\pi_E=q\id_E$,
traduisant le fait que la norme de~$\pi_E$ dans~$\End(E)$ est \'egale \`a~$q$.

\paragraph{Pas de b\'eb\'es, pas de g\'eants}
Pour calculer l'ordre~$d$ d'un \'el\'ement~$g$ d'un groupe ab\'elien fini~$G$,
la m\'ethode la plus \'evidente consiste \`a calculer $g$, $2g$, etc.
jusqu'\`a $dg$, \'egal \`a l'\'el\'ement neutre. Cela requiert $d$ op\'erations
dans le groupe~$G$.

Voici comment D.~\textsc{Shanks}~\cite{shanks1971} propose de proc\'eder 
si l'on conna\^{\i}t un encadrement $D\leq d <D+C$ de~$d$, o\`u $D$ et~$C$
sont des nombres entiers.
Soit $t=\lceil \sqrt C\rceil$ 
le plus petit entier sup\'erieur ou \'egal \`a~$\sqrt C$.
Il existe des entiers~$i$ et~$j$ v\'erifiant $0\leq i,j < t$
tels que $d-D=it+j$, d'o\`u l'\'egalit\'e $(d-D-j)g = i tg$ dans~$G$.

Il suffit alors de calculer d'une part
les~$t$ multiples $g_i=itg$, pour $0\leq i<t$,
de~$tg$, d'autre part les $t$ multiples $b_j=(d-D)g-jg$,
et de d\'eterminer un \'el\'ement de la premi\`ere liste qui appartient
\`a la seconde: si $g_i=b_j$, on a $d=D+it+j$.
Le nom de la m\'ethode, \emph{baby steps---giant steps},
vient de ce que les \'el\'ements~$b_j$ sont dans une progression
de {\og pas de b\'eb\'e\fg}~$g$, tandis que les \'el\'ements~$g_i$
sont dans une progression de {\og pas de g\'eant\fg}~$tg$,
$t$ \'etant approximativement \'egal \`a la racine carr\'ee de~$C$,
suppos\'e grand.
Il en r\'esulte un algorithme 
pour d\'eterminer l'ordre d'un \'el\'ement d'un groupe~$G$,
requ\'erant $\mathrm O(\sqrt C)$ op\'erations dans~$G$ et le
stockage d'autant d'\'el\'ements de~$G$
lorsque l'on sait que cet ordre appartient \`a un intervalle de
longueur~$C$.

Appliqu\'e au groupe~$E(\F)$, on peut ainsi trouver
l'ordre d'un \'el\'ement donn\'e en un temps proportionnel \`a
un multiple de~$\sqrt q$ (multipli\'e par un facteur logarithmique en~$q$,
correspondant \`a la complexit\'e du calcul dans un corps de cardinal~$q$).
Cela ne fournit cependant pas le cardinal de~$E(\F)$.
Toutefois, si $P$ est un point de~$E(\F)$
d'ordre $d$, on a $t_E\equiv q+1\pmod d$ d'apr\`es le th\'eor\`eme de Lagrange,
d'o\`u la valeur exacte de~$t_E$ si $d>4\sqrt q$.

Un tel point n'existe pas toujours. Toutefois, un lemme
de~\textsc{Mestre} affirme que si~$E$ n'a pas
de point d'ordre au moins~$4\sqrt q$, 
alors sa {\og tordue quadratique\fg} $E'$, en poss\`ede un,
tout au moins si $q$ est assez grand ($q\geq 1373$ suffit certainement,
\cf\cite{schoof1995}, th.~3.1).
Si, par exemple, $E$ est donn\'ee par l'\'equation~\eqref{eq.weierstrass},
cette courbe~$E'$ est donn\'ee par l'\'equation
\[ y^2=x^3 + a\omega^2 x+b\omega^3, \]
o\`u $\omega$ est un \'el\'ement de~$\F^*$ qui n'est pas un carr\'e.
(Le moyen le plus simple d'obtenir un tel \'el\'ement~$\omega$
consiste \`a choisir des \'el\'ements de~$\F^*$ au hasard,
jusqu'\`a ce que l'un convienne.)
Les cardinaux de ces courbes sont reli\'es par la relation
\[ \Card {E(\F)}+\Card {E'(\F)} = 2 (q+1), \]
si bien qu'il suffit de d\'eterminer le cardinal de l'une d'entre elles.
En choisissant des points au hasard sur ces deux courbes,
on obtient rapidement un point d'ordre assez grand,
puis, par la m\'ethode de \textsc{Shanks}, le cardinal de~$E(\F)$ et~$E'(\F)$.

Le nombre d'op\'erations que peut requ\'erir l'algorithme de \textsc{Shanks}
est $\mathrm O(q^{1/4})$ et le calcul effectif dans $E(\F)$
n\'ecessite au plus une puissance de~$\log q$ op\'erations \'el\'ementaires.

On a ainsi d\'ecrit un algorithme probabiliste, de complexit\'e 
$\Otilde(q^{1/4})$ pour calculer le nombre de points d'une courbe elliptique sur un corps fini
\`a~$q$ \'el\'ements.

\paragraph{Courbes \`a multiplication complexe}
Supposons que l'on sache que la courbe~$E$ admette
de la multiplication complexe par un ordre~$\mathscr O$
de discriminant~$-D$ d'un corps quadratique imaginaire~$K$.
Comme nous l'avons dit dans le paragraphe consacr\'e \`a l'extraction
de racines carr\'ees, l'endomorphisme de Frobenius est de
norme~$q$ dans~$\End(E)$, donc s'\'ecrit
$\pi_E=\frac12(x+y\sqrt{-D})$ dans~$K$, o\`u $x$, $y$ sont des entiers
relatifs de m\^eme parit\'e qui v\'erifient la relation
\begin{equation}\label{eq.cornacchia}
x^2+Dy^2=4q.
\end{equation}
Alors, $t_E=x$.
Comme $\pi_E$ n'est pas un multiple dans~$\mathscr O$,
le pgcd de~$x$ et~$y$ est~$1$ ou~$2$. Supposons-les impairs
pour simplifier, les adapations \`a faire dans la suite sont \'evidentes.

L'algorithme de \textsc{Cornacchia}
(voir~\cite{nitaj1995} pour une d\'emonstration \'el\'ementaire)
et, en fait, tout algorithme de r\'eduction de r\'eseau fournissant
un vecteur de petite longueur dans un r\'eseau euclidien de dimension~$2$,
permet de calculer tous les couples~$(x,y)$
d'entiers premiers entre eux qui v\'erifient l'\'equation~\eqref{eq.cornacchia}.

La premi\`ere \'etape consiste \`a d\'eterminer les entiers $u$ tels que
$u^2\equiv -D\pmod{4q}$
et $2q<u<4q$.
Pour trouver un tel~$u$, on commence g\'en\'eralement par r\'esoudre
la congruence modulo~$p$ par un algorithme probabiliste,
puis on utilise la m\'ethode de Newton $p$-adique.

Pour chacun de ces entiers~$u$,
appliquons alors l'algorithme d'Euclide au couple $(4q,u)$
et arr\^etons-nous au premier reste~$x$ strictement inf\'erieur \`a~$2\sqrt q$; 
si \mbox{$(4q-x^2)/D$} est le carr\'e d'un nombre entier~$y$, alors
$(\pm x,\pm y)$ est solution de~\eqref{eq.cornacchia}.
Sauf si $D=-1$ o\`u il faut aussi consid\'erer $(\pm y,\pm x)$,
on obtient ainsi toutes les solutions (primitives)
de l'\'equation~\eqref{eq.cornacchia}.

Par cette m\'ethode, on peut donc calculer $\Card{E(\F)}$
au choix du signe de~$x$ (voire celui de~$y$ si $D=-1$) pr\`es.
In\'evitable puisqu'on ne distingue pas
la courbe~$E$ de sa tordue quadratique,
cette ambigu\"{\i}t\'e se r\'esout toutefois sans peine, par exemple
en regardant l'ordre d'un point de~$E(\F)$ pris au hasard.

\subsection{Algorithme de \textsc{Schoof}}

\paragraph{Repr\'esentations}
Fixons une cl\^oture alg\'ebrique $\bar\F$ de~$\F$.
Soit $n$ un entier naturel premier \`a la caract\'eristique de~$\F$
et notons
 $E[ n]$ le sous-groupe de~$E(\bar\F)$ form\'e des points~$P$ 
tels que $n P=0$ ; il est isomorphe \`a~$(\Z/n\Z)^2$.
Tout endomorphisme $u$ de~$E$ laisse stable~$E[n]$,
d'o\`u une action de~$\End(E)$ sur $E[n]$ qui,
si l'on choisit une $(\Z/n\Z)$-base de~$E[n]$,
s'identifie \`a un homomorphisme d'anneaux
\begin{equation}
  \rho_n \colon \End(E) \ra \mathrm M_2(\Z/n\Z).
\end{equation}
En particulier, il correspond \`a $\pi_E$ une 
(classe de conjugaison de) matrice
dont le polyn\^ome caract\'eristique 
est pr\'ecis\'ement $X^2-t_E X+q$, modulo~$n$,
cf. l'\'equation~\eqref{eq.frobenius}.

\paragraph{Principe de l'algorithme de \textsc{Schoof}}
En 1985, R.~\textsc{Schoof}~\cite{schoof1985}
propose de calculer $\Card{E(\F)}$
de la fa\c{c}on suivante:
\begin{enumerate}
\item calculer $t_E$ modulo~$\ell$, pour des nombres premiers
$\ell_1,\ell_2,\dots$ distincts de la caract\'eristique de~$\F$;
\item en d\'eduire, par le th\'eor\`eme chinois,
$t_E$ modulo le produit~$L$ des~$\ell_i$;
\item en d\'eduire $t_E$ gr\^ace \`a l'in\'egalit\'e de Hasse
si $L>4\sqrt q$.
\end{enumerate}

La seconde \'etape est relativement \'evidente, de m\^eme que la derni\`ere
puisqu'il n'y a qu'un seul entier congru \`a~$t_E$ dans l'intervalle
$[q+1-\frac12L,q+1+\frac12L]$ si $L>4\sqrt q$.
Expliquons donc comment calculer $t_E$ modulo~$\ell$
si $\ell$ est un nombre premier.
R\'ecrivons la relation~\eqref{eq.frobenius} en la sp\'ecialisant
aux points de~$E[\ell]$: pour tout $P\in E[\ell]$, on a
\[ \pi_E^2(P) + q P = t_E \pi_E(P). \]
Si l'on trouve un entier $t$ tel que l'on ait $\pi_E^2(P)+qP=t\pi_E(P)$
pour un point  $P\in E[\ell]$, il vient alors $(t-t_E)\pi_E(P)=O$,
d'o\`u $t=t_E\pmod \ell$ si $P$ est d'ordre~$\ell$ (c'est-\`a-dire $P\neq O$).

\paragraph{Polyn\^omes de division}
La multiplication par~$\ell$ dans~$E$ est donn\'ee par
une transformation rationnelle de~$\P^2$, de degr\'e~$\ell^2$,
de la forme $(x:y:z)\mapsto (P_\ell(x,y,z):Q_\ell(x,y,z):R_\ell(x,y,z))$
o\`u $P_\ell$, $Q_\ell$, $R_\ell$ sont
trois polyn\^omes homog\`enes de degr\'es~$\ell^2$ \`a coefficients dans~$\F$,
premiers entre eux.
L'origine~$O$ de~$E$ \'etant l'unique point \`a l'infini de~$E$,
un point $(x:y:z)\in E(\bar\F)$ appartient \`a~$E[\ell]$
si et seulement si $R_\ell(x,y,z)=0$.
Si l'on se restreint au compl\'ementaire de l'origine, 
le sch\'ema des points d'ordre~$\ell$, $E^*_\ell=E[\ell]\cap\A^2$,
est donn\'e par les \'equations $R_\ell(x,y,1)=0$
et~\eqref{eq.weierstrass.general}.

Pour effectuer la premi\`ere \'etape, il faudrait donc
\^etre capable de calculer les coordonn\'ees d'un point de~$E^*_\ell$,
c'est-\`a-dire de trouver une solution explicite $(x,y)\in\bar\F^2$
de ce syst\`eme d'\'equations. Malheureusement, m\^eme s'il existe
des algorithmes de complexit\'e polynomiale en~$\log q$ pour cela,
ces algorithmes sont tous probabilistes, c'est-\`a-dire
qu'ils requi\`erent l'utilisation d'un g\'en\'erateur de nombres al\'eatoires
et que l'on a seulement une tr\`es forte probabilit\'e que l'algorithme
se termine en un temps polynomial en~$\log q$.
\textsc{Schoof} contourne cette difficult\'e en s'int\'eressant directement
\`a la totalit\'e des points d'ordre~$\ell$, c'est-\`a-dire
au sch\'ema~$E^*_\ell$ lui-m\^eme.
Soit $A_\ell$ l'anneau des fonctions de ce sch\'ema~$E^*_\ell$,
c'est-\`a-dire
\begin{equation}
  A_\ell= \F[x,y] / ( y^2+a_1xy+a_3y - x^3-a_2x^2-a_4 x-a_6, R_\ell(x,y,1)).
\end{equation}
L'image du couple d'ind\'etermin\'ees~$(x,y)$  dans~$A_\ell^2$
donne les coordonn\'ees d'un point, tautologique, 
d'ordre~$\ell$ de~$E$ \`a coefficients dans~$A_\ell$.
Les endomorphismes $\pi_E^2+q$ et $t\pi_E$, 
pour $t\in\Z$ non multiple de~$\ell$, 
induisent des endomorphismes de~$E^*_\ell$,
donc des endomorphismes de~$A_\ell$. Vu la pr\'esentation
donn\'ee de l'alg\`ebre~$A_\ell$, 
un endomorphisme $u$ de~$A_\ell$ d\'etermine
deux \'el\'ements $X_u$ et $Y_u$ de~$A_\ell$, correspondant
aux coordonn\'ees de~$u(P)$, lorsque $P$ est
le point tautologique de~$E^*_\ell(A_\ell)$.

\textsc{Schoof} calcule ces \'el\'ements~$(X_u,Y_u)$
lorsque  $u$ est l'un des endomorphismes $v=\pi_E^2+q\id_E$
et $u_t=t\pi_E$, pour $t\in\{\pm 1,\dots,\pm\frac{\ell-1}2\}$.
S'il y a une co\"{\i}ncidence  $(X_v,Y_v)=(X_{u_t},Y_{u_t})$,
alors $t_E\equiv t\pmod\ell$, sinon, $t_E\equiv 0\pmod\ell$.

\paragraph{Complexit\'e}
Il correspond \`a l'endomorphisme~$\pi_E$
un couple $(X_\pi,Y_\pi)$ de~$A_\ell$, image de~$(x^q,y^q)$.
Son calcul revient \`a une \'el\'evation \`a la puissance~$q$
qu'on effectue, par \'el\'evations successives au carr\'e,
en $\mathrm O(\log q)$ multiplications dans~$A_\ell$.
En utilisant les formules d'addition dans~$E$,  
chacun des~$\ell-1$ couples~$(X_{u_t},Y_{u_t})$
\`a consid\'erer demande~$\mathrm O(1)$ multiplications suppl\'ementaires.
Il faut donc effectuer
$\mathrm O(\ell+\log q)$~multiplications dans~$A_\ell$
pour calculer~$t_E\pmod\ell$ selon la m\'ethode d\'ecrite.

L'alg\`ebre~$A_\ell$ est de dimension~$\ell^2-1$ sur~$\F$
et sa pr\'esentation  est assez pratique. Par exemple, si $p>3$
et que~$E$ est donn\'ee sous la forme~\eqref{eq.weierstrass},
elle est de la forme
\begin{equation}
\F[x,y]/(\psi_\ell(x), y^2-x^3-ax-b),
\end{equation}
o\`u $\psi_\ell$ est un polyn\^ome de degr\'e $(\ell^2-1)/2$,
appel\'e \emph{polyn\^ome de division} et tel que,
si $P$ est un point de~$E$ d'abscisse~$x$,
$\psi_\ell(x)$ soit l'abscisse du point~$\ell P$.

Notons $M_q(n)$ le nombre de multiplications dans~$\F$
requises par une multiplication dans une telle alg\`ebre
de dimension~$n$. Na\"{\i}vement, $M_q(n)=\mathrm O(n^2)$
lorsque $n$ tend vers l'infini,
mais la d\'ecouverte de m\'ethodes de multiplication rapide
(\textsc{Karatsuba}, utilisation
de la transform\'ee de Fourier rapide par~\textsc{Sch\"onhage}
et~\textsc{Strassen}) 
a permis de voir
que $M_q(n)=\mathrm O(n^{1+\eps})$ pour tout~$\eps>0$,
ce qu'on notera ici $M_q(n)=\Otilde(n)$.
De m\^eme, si le corps~$\F$ est pr\'esent\'e sous la forme 
$\F_p[t]/(f)$, o\`u $f$ est un polyn\^ome irr\'eductible 
de degr\'e~$d$ \`a coefficients dans~$\F_p$,
la multiplication dans~$\F$ n\'ecessite $\Otilde(d)$ multiplications
dans~$\F_p$.
Finalement, une multiplication dans~$A_\ell$
requiert $\Otilde (\ell^2\log q)$ multiplications \'el\'ementaires.

Ce sont bien s\^ur des \'evaluations asymptotiques
et les constantes implicites dans ces expressions~$\mathrm O$
et~$\Otilde$ sont grandes; 
ainsi, pendant longtemps, les m\'ethodes rapides n'ont \'et\'e
comp\'etitives que pour de grandes valeurs de~$n$.
Apparemment, l'\'evolution r\'ecente des ordinateurs les rend praticables.

En d\'efinitive, le calcul de $t_E$ modulo~$\ell$
a une complexit\'e major\'ee par~$\Otilde(\ell^{3}(\log q)^2)$.

Par ailleurs, le th\'eor\`eme des nombres premiers entra\^{\i}ne
l'existence d'un nombre r\'eel~$c>0$ tel que l'on ait pour tout
nombre r\'eel~$x>2$ la minoration
\footnote{Il semble qu'on puisse prendre $c=2$ si l'on impose
en outre $x>11$.}
%
\begin{equation}
\prod_{\ell<x} \ell \geq  e^{x/c} .
\end{equation}
Par suite, le produit $L_N=\ell_1\dots \ell_N$ des $N$
premiers nombres premiers est sup\'erieur \`a~$4\sqrt q$ si 
$N$ est minor\'e par un multiple de~$\log q$.
Cela fournit finalement  un algorithme d\'eterministe de calcul de~$t_E$,
et donc de~$\Card{E(\F)}$, dont la complexit\'e est
$\Otilde((\log q)^{5})$.

\paragraph{Am\'eliorations d'Atkin et Elkies}\label{subsec.sea}

Comme aime \`a le pr\'esenter R.~\textsc{Schoof}, 
ces am\'eliorations visent \`a utiliser l'action de Frobenius
sur des objets de taille plus petite que le groupe~$E[\ell]$
lui-m\^eme.
O.~\textsc{Atkin}~\cite{atkin1988,atkin1992}
consid\`ere ainsi la droite projective
$(E[\ell]\setminus\{0\})/\F_\ell^*$ \emph{quotient} 
de~$E[\ell]\setminus\{0\}$,
tandis que N.~\textsc{Elkies}~\cite{elkies1991}
utilise (quand ils existent) 
les sous-espaces propres de l'action de Frobenius,
c'est-\`a-dire les sous-sch\'emas en groupes de~$E[\ell]$.

\'Etudions donc plus en d\'etail l'action de~$\pi_E$ sur le groupe~$E[\ell]$
des points de~$\ell$-torsion. Celle-ci s'interpr\`ete comme
un endomorphisme~$\pi_\ell=\rho_\ell(\pi_E)$ 
de~$E[\ell]$, vu comme~$\F_\ell$-espace 
vectoriel de dimension~$2$. Les espaces propres de~$\pi_\ell$
correspondent aux sous-groupes cycliques~$C$ de~$E[\ell]$
qui sont d\'efinis sur le corps~$\F$. \`A un tel sous-groupe
cyclique~$C$ correspond une courbe elliptique $E'=E/C$, d\'efinie sur~$\F$,
et li\'ee \`a~$E$ par une isog\'enie de degr\'e~$\ell$.
Les invariants~$j$ et~$j'$ des courbes~$E$ et~$E'$ fournissent
alors un point $\F$-rationnel~$(j,j')$ de la courbe modulaire~$X_0(\ell)$,
identifi\'ee abusivement \`a 
son image dans~$\P^1\times\P^1=X_0(1)\times X_0(1)$.

Inversement, si $E$ n'est pas supersinguli\`ere et si son invariant~$j$
n'est ni \'egal \`a~$0$ ni \'egal \`a~$1728$
(courbes qu'on qualifiera d'exceptionnelles), \textsc{Atkin} d\'emontre
qu'il correspond \`a tout invariant~$j'\in\F$ tel que $(j,j')\in X_0(\ell)$ 
une courbe elliptique $E'$ sur~$\F$ et une isog\'enie de degr\'e~$\ell$,
$\phi\colon E\ra E'$. Cela repose sur le fait que sur
les $\bar\F$-endomorphismes de telles courbes elliptiques sont d\'efinies
sur~$\F$ et que seuls~$\pm\id$  sont d'ordre fini,
voir~\cite{schoof1995}.

Les courbes exceptionnelles sont trait\'ees ind\'ependamment.
Celles d'invariants~$0$ et~$1728$ admettent
des multiplications complexes par~$\Z[j]$ et~$\Z[i]$
respectivement et leur cardinal se calcule facilement
par l'algorithme de \textsc{Cornacchia}.
Quant aux courbes supersinguli\`eres, leur nombre de points dans~$\F$
est si particulier (\cf\cite{waterhouse1969}, th.~4.1)
que l'on peut rapidement, en choisissant des points
au hasard, v\'erifier si $E$ est supersinguli\`ere et calculer $E(\F)$.
Nous supposons donc que $E$ n'est pas exceptionnelle.

L'\'equation de la courbe~$X_0(\ell)$ dans~$\P^1\times\P^1$
est un polyn\^ome 
$\Phi_\ell\in\Z[X,Y]$, sym\'etrique et de degr\'e~$\ell+1$ 
en chacune des variables ---
le \emph{polyn\^ome modulaire}; autrement dit, $\Phi_\ell(j,X)$
est le polyn\^ome minimal sur~$\Q(j)$ de la fonction m\'eromorphe $j(q^\ell)$
sur $X_0(\ell)$.
Son calcul effectif est possible, au moins si $\ell$ n'est
pas trop grand,
soit \`a l'aide du d\'eveloppement
en s\'erie de Fourier de la fonction modulaire~$j(q)$ et des fonctions~$j(q^\ell)$,
$j(\zeta q)$, pour $\zeta^\ell=1$ (\cf\cite{morain1995}), 
soit par interpolation en choisissant des valeurs particuli\`eres
de~$q$ (\cf\cite{enge2006}).
Toutefois, les polyn\^omes modulaires sont de hauteur tr\`es grande:
d'apr\`es P.~\textsc{Cohen}~\cite{cohen84}, lorsque $\ell$ tend vers l'infini,
\begin{equation}
h(\Phi_\ell) \sim 6 \ell \log \ell\ ;
\end{equation}
pour donner un exemple, 
le terme constant de~$\Phi_5$ poss\`ede~43 chiffres d\'ecimaux!
Comme $\ell$ est de l'ordre de~$\log q$,
on utilise d'autres \'equations de la courbe modulaire~$X_0(\ell)$, 
donn\'ees par le polyn\^ome minimal d'autre fonctions sur~$X_0(\ell)$. 
\textsc{Atkin} a par exemple propos\'e d'employer la fonction
\begin{equation}
f_\ell(\tau) = \ell^s \left( \frac{\eta(\ell\tau)}{\eta(\tau)}\right)^{2s}, 
\qquad s = 12/\pgcd(12,\ell-1), \end{equation}
o\`u $\eta$ d\'esigne la fonction~$\eta$ de \textsc{Dedekind},
donn\'ee par
\begin{equation}
\eta(\tau) = q^{1/24} \prod_{n=1}^\infty (1-q^n) , \qquad q=e^{2i\pi\tau}.
\end{equation}

Si $j$ est l'invariant de la courbe~$E$, suppos\'ee non exceptionnelle,
la factorisation de~$\Phi_\ell(j,Y)$ dans~$\F[Y]$ refl\`ete
donc l'action de~$\pi_E$ sur~$E[\ell]$:
\begin{enumerate}
\item
si $\pi_\ell$ est diagonalisable, $\Phi_\ell(j,Y)$ poss\`ede deux facteurs
de degr\'e~$1$, ses autres facteurs sont de m\^eme degr\'e~$r$,
et $t_2^2-4q$ est un carr\'e modulo~$\ell$;
\item
si $\pi_\ell$ n'est pas semi-simple, $t_E^2-4q\equiv 0\pmod\ell$
et $\Phi_\ell(j,Y)$ poss\`ede exactement deux facteurs irr\'eductibles,
l'un de degr\'e~$1$ et l'autre de degr\'e~$r=\ell$;
\item
si $\pi_\ell$ est semi-simple, non diagonalisable,
$t_E^2-4q$ n'est pas un carr\'e modulo~$\ell$
et les facteurs irr\'eductibles de~$\Phi_j(j,Y)$ sont tous
de m\^eme degr\'e~$r$.
\end{enumerate}
Avec ces notations, $r$ est l'ordre de la matrice~$\pi_\ell$
et il existe un \'el\'ement $\xi\in\bar{\F_\ell}^*$ d'ordre~$r$
tel que $t_E\equiv q(\xi+\xi^{-1})^2\pmod \ell$.

Dans les deux premiers cas, il existe un
sous-groupe~$C$ de~$E[\ell]$ d\'efini sur~$\F$, correspondant
\`a un quotient~$A'_\ell$ de l'alg\`ebre~$A_\ell$ utilis\'ee par~\textsc{Schoof},
c'est-\`a-dire \`a un facteur $\psi_C$ 
de degr\'e~$\ell-1$ du polyn\^ome de division~$\psi_\ell$.
Pour d\'eterminer ce facteur sans expliciter~$\psi_\ell$,
N.~\textsc{Elkies} explique dans~\cite{elkies1991} comment construire
une courbe elliptique~$E'$ sur~$\F$ et une isog\'enie $\phi\colon E\ra E'$
de noyau~$C$.  Ses formules utilisent la th\'eorie des fonctions
elliptiques et font intervenir des d\'enominateurs; elles ne conviennent
que si la caract\'eristique du corps~$\F$ est sup\'erieure \`a~$\ell$,
\cf\cite{schoof1995} et le chapitre~17 de~\cite{cohen2006}.
En {\og petite caract\'eristique\fg}, diverses m\'ethodes existent:
utilisation de la loi de groupe formel ou
du sous-groupe de~$p$-torsion (Couveignes, \cite{lercier-morain2000,couveignes1996}),
\cite{lercier1996} en caract\'eristique~$2$.

Dans le troisi\`eme cas, la factorisation de~$\Phi_\ell$
ne fournit qu'une information partielle sur~$t_E\pmod\ell$,
dont O.~\textsc{Atkin} a montr\'e  comment  la connaissance pouvait 
acc\'el\'erer grandement le calcul de~$t_E$.

\paragraph{R\'esultats}
L'algorithme obtenu en combinant la m\'ethode de~\textsc{Schoof}
et les am\'eliorations d'\textsc{Atkin} et~\textsc{Elkies}
est surnomm\'e~\textsc{sea}. Sa complexit\'e, en temps et en espace,
est polynomiale en~$\log q$,
respectivement $\Otilde((\log q)^{4})$ et~$\mathrm O((\log q)^2)$. 
Toutefois, comme il requiert la factorisation
d'un polyn\^ome \`a coefficients dans un corps fini,
c'est un algorithme probabiliste.
Il est d\'ecrit en grand d\'etail dans~\cite{mueller1995}
lorsque la caract\'eristique du corps n'est pas trop petite;
je renvoie aussi \`a la pr\'esentation de~\textsc{Morain}~\cite{morain1995}.

Son impl\'ementation concr\`ete a \'et\'e r\'ealis\'ee par de nombreuses
personnes, dans de nombreux syst\`emes de calcul alg\'ebrique,
dont Magma~\cite{bosma-c-p1997} et Pari/GP~\cite{pari/gp}.
L'impl\'ementation en caract\'eristique~$2$ de \textsc{Vercauteren}
lui a permis de  calculer le cardinal d'une courbe
elliptique sur~$\F_{2^{1999}}$ \`a l'aide de 10~Pentium~II 400~Mhz 
en environ une semaine,
\cf \cite{vercauteren2000}.
Le record actuel est d\'etenu
par \textsc{Enge}, \textsc{Gaudry} et \textsc{Morain} qui ont calcul\'e
le nombre de points d'une courbe elliptique sur un corps fini
de cardinal le nombre premier~$p=10^{2099}+6243$,
\cf\cite{morain2006}.
Le temps de calcul sur une machine puissante 
(Processeur AMD 64 3400+, $2{,}4$~GHz) est de l'ordre de 200~jours, 
non compris le calcul des polyn\^omes modulaires!

\vskip 0pt plus .1\textheight
\penalty -30
\vskip 0pt plus -.1\textheight

\subsection{G\'en\'eralisations}

\paragraph{Cohomologie \'etale}\label{subsec.l-adique}
Pour expliquer le titre de ce chapitre et le principe des g\'en\'eralisations
de l'algorithme de~\textsc{Schoof}, il nous faut faire quelques rappels
sur la cohomologie~$\ell$-adique.
Soit $\F$ un corps fini, notons $q$ son cardinal. 
Fixons une cl\^oture alg\'ebrique~$\bar\F$ de~$\F$ et notons
$\F_{q^n}$ le sous-corps \`a~$q^n$ \'el\'ements de~$\bar\F$, de sorte que $\F_q=\F$.
Soit $\phi$ l'automorphisme
de Frobenius g\'eom\'etrique de~$\bar\F$, inverse de l'automorphisme
de Frobenius arithm\'etique $x\mapsto x^q$.

Soit $\ell$ un nombre premier 
distinct de la caract\'eristique de~$\F$;
notons $\F_\ell$ le corps~$\Z/\ell\Z$,
$\Z_\ell$ l'anneau des entiers~$\ell$-adiques
et $\Q_\ell$ son corps des fractions.

Soit $X$ un sch\'ema s\'epar\'e de type fini sur~$\F$, posons 
$\bar X=X\otimes\bar\F$.
Les groupes de cohomologie $\ell$-adique  \`a support propre de~$\bar X$
d\'efinis par 
\textsc{Grothendieck},
\begin{equation}
 H^i_c (\bar X,\F_\ell), \quad H^i_c(\bar X,\Q_\ell), \quad H^i_c(\bar X,\Z_\ell), \end{equation}
sont respectivement des espaces vectoriels de  dimension finie
sur~$\F_\ell$, $\Q_\ell$, et des $\Z_\ell$-modules de type fini,
nuls si $i<0$ ou si $i>2\dim X$.
Ils sont munis d'une action de~$\Gal(\bar\F/\F)$ et permettent le calcul
du cardinal de~$X(\F)$
via une formule de Lefschetz~\cite{grothendieck1964,sga5} :
\begin{equation}\label{eq.lefschetz.Ql}
\Card{X(\F)} = \sum_{i=0}^{2\dim X} (-1)^i \Tr(\phi| H^i_c(\bar X,\Q_\ell)). 
\end{equation}
Cette formule, appliqu\'ee aux extensions finies de~$\F$, entra\^{\i}ne
la formule suivante pour la fonction z\^eta de~$X$:
\begin{equation}\label{eq.zeta}
 Z(X,t) = \exp\big(\sum_{n=1}^\infty \Card {X(\F_{q^n})} \frac{t^n}n\big)
  =  \prod_{i=0}^{2\dim X} \det(1-t\phi|H^i_c(\bar X,\Q_\ell)^{(-1)^{i+1}},
\end{equation}
d'o\`u, de nouveau, la rationalit\'e de  la fonction z\^eta de~$X$.

On a $H^i_c(\bar X,\Q_\ell)=H^i_c(\bar X,\Z_\ell)\otimes_{\Z_\ell}\Q_\ell$,
et s'il n'est pas vrai que $H^i_c(\bar X,\Z_\ell)\otimes_{\Z_\ell}\F_\ell$
est \'egal \`a $H^i_c(\bar X,\F_\ell)$, ces deux groupes peuvent \^etre
reli\'es, de sorte que l'on a une congruence modulo~$\ell$
(\cite{deligne1977}, Fonctions~$L$ modulo~$\ell^n$, p.~116, th.~2.2) :
\begin{equation}\label{eq.lefschetz.Fl}
\Card{X(\F)} \equiv  \sum_{i=0}^{2\dim X} \Tr(\phi| H^i_c(\bar X,\F_\ell))
 \pmod\ell.
\end{equation}

Supposons par exemple que $X=E$ soit une courbe elliptique~$E$.
Si $A$ est l'un des anneaux $\F_\ell$, $\Z_\ell$, $\Q_\ell$,
alors $H^0_c(\bar E,A)=A$, $H^2_c(\bar E,A)=A$, les endomorphismes~$\phi$
\'etant respectivement l'identit\'e et la multiplication par~$q$,
tandis que $H^1_c(\bar E,A)$ est un $A$-module libre de rang~$2$.
En outre, $H^1_c(\bar E,\F_\ell)$, muni de l'action de~$\phi$,
s'identifie canoniquement \`a~$E[\ell]$, 
muni de l'action de l'endomorphisme~$\pi_E$,
et le polyn\^ome $X^2-t_E X+q$ est le polyn\^ome caract\'eristique
de~$\phi$ agissant sur~$H^1_c(E,\F_\ell)$.
En particulier, $t_E\equiv \Tr(\phi|H^1_c(E,\F_\ell)) \pmod\ell$.
Dans ce cas, l'entier~$t_E$ est d'ailleurs la trace commune
de~$\phi$ sur tous les espaces $H^1_c(E,\Z_\ell)$,
$H^1_c(E,\Q_\ell)$ et l'on a
\begin{equation}
Z(E,t)= \frac{1-t_E t + q t^2}{(1-t)(1-qt)}.
\end{equation}

Ainsi, avec ces identifications, la congruence~\eqref{eq.lefschetz.Fl}
n'est autre que celle qui est \`a la base de l'algorithme de~\textsc{Schoof},
d'o\`u le titre, \emph{approche $\ell$-adique}, de ce chapitre.

Revenons au cas g\'en\'eral en supposant que $X$ soit projective
et lisse.
P.~\textsc{Deligne}~\cite{deligne74} a d\'emontr\'e
que pour tout~$i$,
le polyn\^ome caract\'eristique de~$\phi$ agissant
sur l'espace $H^i_c(\bar X,\Q_\ell)$ est un polyn\^ome
\`a coefficients entiers qui ne d\'epend pas de~$\ell$
(distinct de la caract\'eristique de~$\F$)  
et dont les racines complexes sont toutes de module $q^{i/2}$.
Omettant par abus l'anneau~$\Q_\ell$ des notations,
on a en particulier la majoration
\begin{equation} \label{eq.majoration.weil}
\abs{\Tr(\phi| H^i_c(\bar X))}  \leq  q^{i/2} \dim H^i_c(\bar X) .
\end{equation}
Cette in\'egalit\'e g\'en\'eralise la majoration de~\textsc{Hasse};
conjectur\'ee par \textsc{Weil},
c'est l'analogue pour la vari\'et\'e~$X$ de l'hypoth\`ese de Riemann.

Pour calculer $\Card{X(\F)}$, 
il suffit de calculer $\Tr(\phi|H^i_c(\bar X))$.
Faisons l'hypoth\`ese que $X$
est projective, lisse et g\'eom\'etriquement int\`egre
(sans pour autant supprimer l'indice~$c$ de la cohomologie),
et supposons que l'on a 
$H^i_c(\bar X,\Z_\ell))\otimes\F_\ell \simeq H^i_c(\bar X,\F_\ell)$.
On a alors la congruence
\[ \Tr(\phi|H^i_c(\bar X))\equiv \Tr(\phi|H^i_c(\bar X,\F_\ell))
\pmod\ell  . \]
Adapt\'e \`a ce cadre, l'algorithme de~\textsc{Schoof} 
proc\'ederait de la fa\c{c}on suivante:
\begin{enumerate}
\item calculer $\Tr(\phi|H^i_c(\bar X,\F_\ell))$ pour 
$\ell=\ell_1,\ell_2,\dots,\ell_n$ et $0\leq i\leq 2\dim X$;
\item en d\'eduire, par le th\'eor\`eme chinois, 
l'entier~$\Tr(\phi|H^i_c(\bar X))$
modulo $L=\ell_1\dots\ell_n$;
\item en d\'eduire l'entier~$\Tr(\phi|H^i_c(\bar X))$
si $L$ est au moins \'egal \`a deux fois la borne
donn\'ee par l'\'equation~\eqref{eq.majoration.weil}.
\end{enumerate}

Le probl\`eme est d'avoir une prise raisonnable sur le
groupe de cohomologie $H^i_c(\bar X,\F_\ell)$.
Sous les hypoth\`eses donn\'ees,
la dualit\'e de Poincar\'e identifie $H^i_c$ muni de~$\phi$
et $H^{2\dim X-i}_c$ muni de $q^{\dim X}\phi^{-1}$.
Supposons donc $0\leq i\leq \dim X$.
Pour $i=0$, on a $H^0_c(\bar X,\F_\ell)=\F_\ell$  et $\phi=\id$;
la trace cherch\'ee vaut~$1$.
Pour $i=1$, $H^1_c(\bar X,\F_\ell)$ s'identifie au groupe correspondant
$H^1_c(A,\F_\ell)$, o\`u $A$ est la vari\'et\'e d'Albanese de~$X$.
C'est une vari\'et\'e ab\'elienne que l'on peut esp\'erer d\'ecrire
explicitement de sorte \`a appliquer la m\'ethode de~\textsc{Schoof}.
Toutefois, pour $2\leq i\leq \dim X$, il ne semble pas y avoir,
en g\'en\'eral, de description raisonnablement effective
de $H^i_c(\bar X,\F_\ell)$. Cette approche restera donc limit\'ee
aux courbes et aux vari\'et\'es ab\'eliennes, dont la cohomologie
est contr\^ol\'ee par le~$H^1$,
voire aux vari\'et\'es pour lesquelles l'on peut d\'ecrire 
effectivement et efficacement
la cohomologie \`a l'aide de vari\'et\'es ab\'eliennes.

\paragraph{Courbes et vari\'et\'es ab\'eliennes}
G\'en\'eralisant des r\'esultats de J.~\textsc{Pila}~\cite{pila1990},
L.~\textsc{Adleman} et M.-D.~\textsc{Huang}
ont ainsi donn\'e dans~\cite{adleman-huang2001} un algorithme
qui calcule le nombre de points d'une vari\'et\'e ab\'elienne~$A$ de dimension~$g$
d\'efinie sur un corps fini~$\F_q$ \`a~$q$ \'el\'ements, plong\'ee dans l'espace
projectif~$\P^N$. Il convient de remarquer
que ces algorithmes ne calculent pas les traces $\Tr(\phi|H^i_c(\bar A,\F_\ell))$
individuellement, mais leur somme altern\'ee, c'est-\`a-dire 
une congruence $\card{A(\F_q)}\pmod\ell$.
Supposant que l'on dispose
de polyn\^omes homog\`enes $F_1,\dots,F_S$ de degr\'es au plus~$T$
d\'efinissant l'id\'eal de~$A$ dans~$\F_q[x_0,\dots,x_N]$,
de formules pour l'addition de~$A$ donn\'ees
par des polyn\^omes de degr\'e~$D$ dans~$R$ cartes affines,
la complexit\'e de cet algorithme,
$(\log q)^{\mathrm O(N^2(g+\log R+\log D))}$.
L'exposant de~$\log q$ est \'enorme : si, pour simplifier,
le plongement est donn\'e par le cube d'une polarisation principale,
on a $N+1=3^g g!$. 


Par ces techniques, il est aussi possible de calculer le polyn\^ome
caract\'eristique de l'endomorphisme de Frobenius~$\pi_A$
agissant sur le module de Tate~$T_\ell(A)$ de~$A$ (pour $\ell$
distinct de la caract\'eristique de~$\F_q$).
\textsc{Pila} proc\`ede par exemple en calculant, pour
tout polyn\^ome irr\'eductible $r\in(\Z/\ell\Z)[t]$, de degr\'e au plus~$2\dim A$,
et leurs puissances,
le cardinal du noyau de l'endomorphisme~$r(\phi)$  du groupe fini~$A[\ell]$.
Il en d\'eduit le polyn\^ome caract\'eristique de~$\phi$ sur le $F_\ell$-espace
vectoriel~$A[\ell]$, puis
le polyn\^ome caract\'eristique de~$\phi$
via le th\'eor\`eme chinois, ayant choisi des valeurs
de~$\ell\leq\mathrm O(\log q)$.

Supposons donc que $A$ soit la jacobienne d'une courbe~$C$ sur~$\F_q$,
suppos\'ee projective, lisse et g\'eom\'etriquement int\`egre.
D'apr\`es \textsc{Weil}~\cite{weil1948}, on a l'\'egalit\'e 
\begin{equation}
 \Card{C(\F_q)} = q + 1 - \Tr(\pi_A). \end{equation}
Comme on peut d\'ecrire effectivement une jacobienne,
il existe donc un algorithme de complexit\'e polynomiale en~$\log q$
pour calculer le nombre de points d'une telle courbe sur un corps fini.
D'apr\`es~\textsc{Pila}, la complexit\'e d'un tel algorithme
est uniforme lorsque la courbe parcourt une famille alg\'ebrique,
voir~\cite{pila1991}.
Si elle reste trop grande pour que cette m\'ethode
puisse \^etre utilis\'ee en pratique,
le cas des courbes hyperelliptiques  a fait l'objet d'am\'eliorations
importantes.

Supposons que $A$ soit la jacobienne d'une courbe hyperelliptique~$C$ qui
est donn\'ee par une \'equation sous la forme (affine)
\begin{equation}\label{eq.hyperelliptique}
y^2 = f(x), \qquad f\in \F_q[x], \quad \deg f=2g+1. 
\end{equation}

Toute classe de diviseur de degr\'e~$0$ sur~$C$ est alors
repr\'esent\'ee de mani\`ere unique par deux polyn\^omes~$u$ et~$v$
v\'erifiant les conditions suivantes:
\begin{enumerate}\def\theenumi{\roman{enumi}}\def\labelenumi{(\theenumi)}
\item $u$ est unitaire;
\item $\deg v<\deg u\leq g$;
\item $u$ divise~$v^2-f$.
\end{enumerate}
L'id\'eal du diviseur dans la carte affine ci-dessus n'est autre
que $(u(x),y-v(x))$.
Cette repr\'esentation,
rappel\'ee dans~\cite{mumford1984}, est 
souvent appel\'ee \emph{description de Mumford} dans la litt\'erature.
Elle permet \`a \textsc{Adleman} et~\textsc{Huang}
de montrer l'existence d'un algorithme pour calculer le cardinal de~$C(\F_q)$
et de~$A(\F_q)$ dont la complexit\'e est $(\log q)^{\mathrm O(g^2\log g)}$.

Bien que consid\'erablement inf\'erieure \`a celle des vari\'et\'es
ab\'eliennes g\'en\'erales, cette complexit\'e  reste exponentielle en le genre,
et ces m\'ethodes sont impropres aux applications cryptographiques.
Voir toutefois l'article~\cite{gaudry-harley2000}
par \textsc{Gaudry} et~\textsc{Harley}
concernant les courbes de genre~$2$ : l'usage de l'analogue des
polyn\^omes de division introduits par D.~\textsc{Cantor}
(voir~\cite{cantor1994})
leur permit de calculer le nombre de points d'une telle courbe
sur un corps fini dont le cardinal est de l'ordre de~$10^{20}$.

\subsection{Formes modulaires}
Revenons pour l'instant au cas des courbes elliptiques
et supposons que $E$ soit une courbe elliptique d\'efinie sur~$\Q$,
donn\'ee par une \'equation de Weierstrass~\eqref{eq.weierstrass},
o\`u $a$ et~$b$ sont des entiers relatifs.
Pour tout nombre premier~$p$, on peut r\'eduire l'\'equation modulo~$p$
et en d\'eduire, tout au moins si $p$ ne divise pas le discriminant
$\Delta_E=-4a^3-27b^2$, une courbe elliptique sur~$\F_p$; \'ecrivons son 
cardinal sous la forme $p+1-a_p$. Pour les quelques nombres premiers qui
restent, on peut d\'efinir un entier~$a_p$ analogue et d\'efinir
la fonction~$L$ de Hasse-Weil de~$E$ par le produit eul\'erien 
et la s\'erie de Dirichlet
\begin{equation}
L(E,s)=\prod_{\substack{\text{$p$ premier} \\ p\nmid \Delta_E}}
                    (1-a_p p^{-s}+p^{1-2s})^{-1}
       \prod_{p\nmid\Delta_E} (1-a_p p^{-s})^{-1}
= \sum_{n=1}^\infty a_n n^{-s}. 
\end{equation}
Ce produit et cette s\'erie convergent pour $\Re(s)>3/2$; d'apr\`es le th\'eor\`eme
de~\textsc{Wiles} et~\textsc{Taylor--Wiles}, compl\'et\'e par~\cite{bcdt2001},
ils poss\`edent un prolongement holomorphe \`a~$\C$ et une \'equation fonctionnelle
reliant $L(E,s)$ \`a~$L(E,2-s)$.
Plus pr\'ecis\'ement, si $N_E$ d\'esigne le conducteur de~$E$,
la fonction holomorphe sur le demi-plan de Poincar\'e
donn\'ee par le d\'eveloppement de Fourier
\begin{equation}\label{eq.fourier}
f_E(\tau)= \sum_{n=1}^\infty a_n q^n, \qquad q=e^{2i\pi\tau} ,
\end{equation}
est une forme modulaire de poids~$2$ pour le sous-groupe de congruence
$\Gamma_0(N_E)$ de~$\SL_2(\Z)$.
L'algorithme de~\textsc{Schoof} appara\^{\i}t ainsi comme un algorithme
de calcul des coefficients des formes modulaires de poids~$2$.

Dans un article r\'ecent~\cite{edixhoven-cjmb2006},
\textsc{Edixhoven}, en collaboration
avec \textsc{Couveignes}, \textsc{de Jong}, \textsc{Merkl} et~\textsc{Bosman},
explique comment calculer les coefficients 
de Fourier d'une forme modulaire de poids~$k\geq 2$, niveau~$N$, parabolique
et propre pour les op\'erateurs de Hecke. Leur approche n'est cependant
enti\`erement men\'ee au bout que pour la fonction~$\Delta$ de Ramanujan,
donn\'ee par
\begin{equation}\label{eq.Delta}
\Delta(\tau) = \eta(\tau)^{24} = q \prod_{n=1}^{\infty} (1-q^n)^{24}
= \sum_{n=1}^\infty a_n q^n,
\qquad q=e^{2i\pi\tau}.
\end{equation}
(L'entier~$a_n$ est classiquement not\'e~$\tau(n)$.)
Je me limite ici \`a une description rapide de quelques-unes
des id\'ees essentielles de ce long article. Soit donc $f$ une 
forme modulaire sur~$\Gamma_1(N)$, propre pour les op\'erateurs de Hecke,
normalis\'ee, de d\'eveloppement de Fourier $\sum a_n q^n$.

\begin{enumerate}
\item Si $f$ \'etait une s\'erie d'Eisenstein, son coefficient~$a_n$
serait une somme explicite de puissances de diviseurs de~$n$; on ne sait pas
\'evaluer une telle somme sans factoriser~$n$, et l'on ne sait
pas factoriser~$n$ en temps polynomial en~$\log n$.
On se contentera donc des coefficients~$a_p$, pour $p$ un nombre premier.

\item D'apr\`es~\textsc{Deligne}, voir~\cite{deligne71b},
il existe pour tout~$\ell$
une repr\'esentation de degr\'e~$2$ de~$\Gal(\bar\Q/\Q)$,
\`a coefficients dans~$\Z_\ell$,
\begin{equation}
\rho_{f,\ell}\colon\Gal(\bar\Q/\Q)\ra\GL_2(\Z_\ell)
\end{equation}
telle que pour tout~$p\neq\ell N$,
$a_p$ soit la trace de~$\rho_f(\phi_p)$, o\`u $\phi_p$
est un \'el\'ement de Frobenius (g\'eom\'etrique) en la place~$p$.
La construction g\'eom\'etrique de cette repr\'esentation et la d\'emonstration
par~\textsc{Deligne} des conjectures de Weil entra\^{\i}nent en outre la majoration
\begin{equation}
\abs{a_p}  \leq 2 p^{(k-1)/2}
\end{equation}
qu'avait conjectur\'ee Ramanujan.

La m\'ethode de~\textsc{Schoof} sugg\`ere donc de calculer 
la r\'eduction modulo~$\ell$, disons $\bar\rho_{f,\ell}$,
de cette repr\'esentation pour
des valeurs de~$\ell$ au plus \'egale \`a~$\mathrm O(\log p)$.

\item
Lorsque $k=2$, cette repr\'esentation~$\bar\rho_{f,\ell}$
se r\'ealise dans celle associ\'ee aux points de~$\ell$-torsion
de la jacobienne~$J_1(N)$ de la courbe modulaire~$X_1(N)$.
Dans le cas g\'en\'eral, la construction de~\textsc{Deligne}
fait intervenir un groupe de cohomologie de degr\'e~$k-1$
d'un produit fibr\'e (d\'esingularis\'e) $k-2$-fois
de la {\og courbe elliptique universelle\fg} sur $X_1(N)$.
Comme on l'a \'evoqu\'e plus haut, il ne semble pas possible de 
d\'ecrire cette cohomologie explicitement.

En revanche, des ph\'enom\`enes de \emph{congruence} entre formes modulaires
entra\^{\i}nent que~$\bar\rho_{f,\ell}$ se r\'ealise dans la repr\'esentation
galoisienne associ\'ee aux points de~$\ell$-torsion
de la jacobienne~$J_1(N\ell)$. 

Comme le genre de~$J_1(N\ell)$ est de l'ordre de~$N^2\ell^2$,
il s'agit de d\'etecter une sous-repr\'esentation de dimension~$2$ dans
une repr\'esentation de tr\`es grande dimension.
Autrement dit, de d\'etecter les $\ell^2-1$ points d'ordre~$\ell$
de~$J_1(N\ell)$ correspondant \`a~$\bar\rho_{f,\ell}$,
le polyn\^ome minimal d'un g\'en\'erateur de l'extension 
de~$\Q$ engendr\'ee par leurs coordonn\'ees et l'action
du groupe de Galois~$\Gal(\bar\Q/\Q)$.

\item
Pour ce faire, \textsc{Couveignes} 
a sugg\'er\'e de travailler dans~$\C$ et de calculer une approximation
de ce polyn\^ome, en m\^eme temps qu'une borne pour sa hauteur. 
(Deux nombres rationnels distincts $x=a/b$ et $x'=a'/b'$ de hauteurs
$\log\max(\abs a,\abs b)$  et $\log\max(\abs{a'},\abs{b'})$
au plus~$h$ diff\`erent d'au moins $e^{-2h}$.)

Les bornes requises sur la hauteur sont obtenues en majorant certaines
quantit\'es concernant  la th\'eorie d'Arakelov des courbes modulaires $X_1(\ell)$:
la hauteur de Faltings, certaines fonctions~$\theta$ et les fonctions
de Green ; ces bornes sont polynomiales en~$\ell$.

\item
Plut\^ot que calculer dans la jacobienne $J_1(N\ell)$, 
les auteurs pr\'ef\`erent utiliser le produit sym\'etrique~$X_1(N\ell)^{(g)}$
de la courbe modulaire, qui est muni d'une application naturelle
birationnelle vers~$J_1(N\ell)$. Pour garantir que 
les points d'ordre~$\ell$ intervenant dans~$\bar\rho_{f,\ell}$
ont un unique ant\'ec\'edent, la construction g\'eom\'etrique d'un
diviseur ayant des propri\'et\'es sp\'ecifiques 
est n\'ecessaire, si bien qu'\`a ce stade, les auteurs de~\cite{edixhoven-cjmb2006}
se cantonnent \`a la forme modulaire~$\Delta$ de Ramanujan.
Pour r\'ealiser cette construction, ils r\'eduisent
modulo~$p$, et l'utilisation d'un algorithme probabiliste leur
est pour l'instant n\'ecessaire.
\end{enumerate}

Tout ceci combin\'e, 
\textsc{Edixhoven} et al. d\'emontrent qu'il existe
un algorithme probabiliste de complexit\'e polynomiale en~$\ell$
calculant $\bar\rho_{f,\ell}$, c'est-\`a-dire :
\begin{enumerate}
\item un corps de nombres $K_\ell$
donn\'e par sa table de multiplication, galoisien sur~$\Q$ correspondant
au noyau de la repr\'esentation~$\bar\rho_{f,\ell}$, 
\item des matrices
correspondant aux \'el\'ements~$\sigma$ du groupe $\Gal(K_\ell/\Q)$
agissant $\Q$-lin\'eairement sur~$K_\ell$,   et
\item 
les matrices $\bar\rho_{f,\ell}(\sigma)\in\GL_2(\F_\ell)$ correspondantes.
\end{enumerate}
Ils en d\'eduisent qu'il existe un algorithme probabiliste
calculant le coefficient~$\tau(p)$ de la fonction de Ramanujan,
algorithme dont la complexit\'e est polynomiale en~$\log p$.

\section{M\'ethodes $p$-adiques}\label{sec.p}

\subsection{Application de la th\'eorie de \textsc{Dwork}}

Comme je l'ai dit plus haut, la premi\`ere d\'emonstration
de la rationalit\'e de la fonction z\^eta d'une vari\'et\'e
alg\'ebrique d\'efinie sur un corps fini est due \`a~\textsc{Dwork}~\cite{dwork1960}.
A.~\textsc{Lauder} et D.~\textsc{Wan}
ont observ\'e~\cite{lauder-wan2006,wan2006}
que la d\'emonstration de~\textsc{Dwork}
permet un moyen de calcul effectif de la fonction z\^eta.
M\^eme si cette m\'ethode s'av\`ere moins efficace que celles qui
ont \'et\'e d\'evelopp\'ees peu apr\`es, elle fournit un algorithme
ind\'ependant de la g\'eom\'etrie du syst\`eme d'\'equations consid\'er\'e.

\paragraph{Calculabilit\'e de la fonction z\^eta}
Soit $\F$ un corps fini, notons $q$ son cardinal et $p$ sa caract\'eristique.
Consid\'erons une vari\'et\'e alg\'ebrique affine~$X$,
lieu des z\'eros dans~$\A^n$ de $m$~polyn\^omes 
$f_1,\dots,f_m\in\F[x_1,\dots,x_n]$ ; soit $a$ le maximum des degr\'es
des~$f_i$.
L'espace n\'ecessaire pour \'ecrire ces polyn\^omes est de l'ordre
de~$m(a+1)^n\log q$.
On cherche un algorithme efficace,
pour calculer la fonction z\^eta de~$X$.

D'apr\`es~\textsc{Dwork},
cette fonction z\^eta s'\'ecrit comme le quotient $P_1/P_2$ de deux polyn\^omes
\`a coefficients entiers, premiers entre eux, de termes constants~$1$.
Un th\'eor\`eme de~\textsc{Bombieri}
majore par~$A=(4a+9)^{n+m}$ la somme des degr\'es de~$P_1$ et de~$P_2$
(prop.~4.2 et th.~1 de~\cite{bombieri1978}, appliqu\'e \`a $f=1$).
Autrement dit, il existe des nombres entiers~$a_1$ et~$a_2$
tels que $a_1+a_2\leq A$, des entiers alg\'ebriques $\alpha_1,\dots,\alpha_{a_1}$,
$\beta_1,\dots,\beta_{a_2}$ tels que pour tout entier~$k\geq 1$,
\begin{equation}\label{eq.Card(X)}
\Card{X(\F_{q^k})} = - \sum_{j=1}^{a_1} \alpha_j^k + \sum_{j=1}^{a_2} \beta_j^k,
\end{equation}
la fonction z\^eta elle-m\^eme \'etant donn\'ee par
\begin{equation} 
Z(X,t) = \frac{ \prod_{j=1}^{a_1} (1-\alpha_jt)}{\prod_{j=1}^{a_2}(1-\beta_jt)}.
\end{equation}
Par cons\'equent, pour calculer la fonction z\^eta de~$X$, il suffit 
de calculer les cardinaux $\Card{X(\F_{q^k})}$ pour $1\leq k\leq A+1$,
d'en d\'eduire les~$A+1$ premiers termes du d\'eveloppement
en s\'erie enti\`ere de~$Z(X,t)$ puis de calculer une approximante
de Pad\'e de cette s\'erie. 
En particulier, \emph{la fonction z\^eta est effectivement calculable.}

D'autre part, d'apr\`es~\textsc{Deligne}~\cite{deligne1980},
les valeurs absolues des~$\alpha_j$ et~$\beta_j$ sont de la forme~$q^{s/2}$,
o\`u $s$ est un entier compris entre~$0$ et~$2\dim X$ (d\'ependant de~$j$).
En particulier, les coefficients de~$P_1$ et de~$P_2$ sont 
major\'es par~$2^Aq^{nA}$.
L'espace n\'ecessaire pour \'ecrire la fonction z\^eta est
donc $\mathrm O(nA\log q)=\mathrm O(n(4a+9)^{n+m}\log q)$, 
sensiblement du m\^eme  ordre de grandeur que celui n\'ecessaire au stockage des
donn\'ees.

La question algorithmique qui se pose naturellement alors est
la suivante:
\emph{Est-il possible de calculer la fonction z\^eta en temps
\emph{polynomial} en cette taille?}
Compte-tenu des bornes de~\textsc{Bombieri} rappel\'ees ci-dessus, 
cela revient \`a la question:
\emph{Est-il possible de calculer le cardinal de~$X(\F_{q^k})$
en temps polynomial en~$k (a+1)^{n+m}\log q$?}

Pour les vari\'et\'es de dimension z\'ero de la droite affine, d\'efinies
par un polyn\^ome $f\in\F[x]$, de degr\'e~$a$, l'algorithme de \textsc{Berlekamp}
a la complexit\'e voulue. Si $f$ est sans racines multiples, le pgcd
de~$f$ et~$x^q-x$ a pour degr\'e le nombres de racines de~$f$
dans~$\F_q$. On calcule bien s\^ur ce pgcd par l'algorithme d'Euclide,
en commen\c{c}ant par \'evaluer $x^q$ dans l'alg\`ebre~$\F[x]/(f)$
ce qui requiert en gros~$\log q$ op\'erations dans cette alg\`ebre; 
la suite de l'algorithme n\'ecessite au plus~$a$ divisions
euclidiennes de polyn\^omes de degr\'es au plus~$a$ \`a coefficients
dans~$\F$.

Pour les courbes planes toutefois, les r\'esultats du chapitre pr\'ec\'edent
faisaient appara\^{\i}tre une d\'ependance exponentielle en le genre...

Si $q=p^d$, \textsc{Lauder} et~\textsc{Wan} montrent 
que la r\'eponse est oui si l'on se limite aux
corps de caract\'eristique major\'ee. Ils d\'eduisent
de l'\'etude par~\textsc{Dwork} de la fonction
z\^eta d'une hypersurface une \emph{congruence} modulo
une puissance de~$p$ pour le nombre
de points d'une hypersurface sur~$\F_q$, d'o\`u sa valeur
si l'exposant de la puissance d\'epasse~$dn$.

\paragraph{Une formule de congruence}
Nous allons nous borner au cas d'une hypersurface de l'espace
affine d'\'equation~$f=0$, o\`u $f\in\F_q[x_1,\dots,x_n]$.
La m\'ethode de~\textsc{Dwork} est {\og torique\fg} 
et commence en fait par \'etudier
la fonction z\^eta de l'intersection~$X$ de cette hypersurface
avec le tore~$\gm^n$.  Le reste s'\'etudie classiquement
par r\'ecurrence, au moyen de la formule d'inclusion-exclusion,
mais cela induit inexorablement
des facteurs~$2^n$ dans la complexit\'e des algorithmes \`a venir.

La th\'eorie des sommes de caract\`eres fournit la formule
\begin{equation}\label{eq.somme.caracteres}
q \Card{X(\F_q)}-(q-1)^n = \sum_{(x_0,\dots,x_n)\in\gm^{n+1}(\F_{q})}
\chi_q(x_0f(x_1,\dots,x_n)), 
\end{equation}
o\`u $\chi_q\colon \F_{q}\ra\Omega^*$ est un caract\`ere
additif non trivial de~$\F_{q}$, $\Omega$ \'etant un corps
de caract\'eristique z\'ero.

 Soit $W$ l'anneau des vecteurs de Witt de~$\F$.
C'est un anneau de valuation discr\`ete, complet, de corps r\'esiduel~$\F$
et dont l'id\'eal maximal est engendr\'e par~$p$.
Si $\F=\F_p$, $W$ n'est autre que l'anneau~$\Z_p$ des entiers~$p$-adiques;
si $\F$ est d\'ecrit sous la forme~$\F_p[x]/(f)$, o\`u $f$ est un 
polyn\^ome de degr\'e~$d$, $W$ est isomorphe \`a l'anneau
quotient~$\Z_p[x]/(\tilde f)$, o\`u $\tilde f\in\Z_p[x]$ est un polyn\^ome
de degr\'e~$d$ arbitraire dont la r\'eduction modulo~$p$ est \'egale \`a~$f$.

Soit $R$ l'anneau obtenu en adjoignant \`a~$W$ 
un \'el\'ement~$\pi$ tel que~$\pi^{p-1}=-p$.
Notons $\Omega$ le corps des fractions de~$R$.
Sa valeur absolue~$\abs\cdot$
et sa valuation~$\ord$ sont normalis\'ees par $\abs p=1/p$
et~$\ord p=1$.  Soit $d$ l'entier tel que $q=p^d$.
Notons $\omega\colon\F\ra\Omega$ 
le caract\`ere de Teichm\"uller,
prolong\'e par~$\omega(0)=0$. Pour tout $x\in\F^*$, $\omega(x)$
est l'unique racine de l'unit\'e d'ordre premier \`a~$p$ de~$W$
qui est congrue \`a~$x$ modulo~$p$; c'est aussi la limite de
la suite $\xi^{q^n}$, o\`u $\xi$ est un \'el\'ement arbitraire de~$W$
qui est congru \`a~$x$ modulo~$p$.

Soit $\sigma$ l'unique automorphisme de~$K$
tel que~$\sigma(\pi)=\pi$ et est congru \`a l'automorphisme
de Frobenius modulo~$p$; notons~$\tau$ l'inverse de~$\sigma$.

Si les coefficients de la s\'erie exponentielle sont $p$-adiquement
trop grands, ceux de la s\'erie~$\theta\in\Omega[[z]]$ d\'efinie par
\begin{equation}\label{eq.theta.dwork}
\theta(z)=\exp(\pi z-\pi z^p)
\end{equation}
sont de valeurs absolues inf\'erieures ou \'egales \`a~$1$
et tendent vers~$0$. Pr\'ecis\'ement,
on a
\begin{equation}\label{eq.major.theta}
\ord \theta_n \geq \max\big( \frac{n(p-1)}{p^2}, \frac2{p-1} \big).
\end{equation}
\textsc{Dwork} montre alors que l'on peut choisir pour caract\`ere~$\chi_q$
dans la formule~\eqref{eq.somme.caracteres} la fonction
\begin{equation}
\chi_q(z) = \theta(1)^{\Tr_{\F_{q}/\F_p}(z)}
 = \prod_{i=0}^{d-1}\theta(\omega(z)^{p^i}).
\end{equation}

Soit $F$ la s\'erie formelle en les ind\'etermin\'ees $x_0,\dots,x_n$
d\'efinie par
\begin{equation}
F(x)=\prod_J \theta(\omega(f_J) x^J), \qquad \text{o\`u}\quad 
 x_0f(x_1,\dots,x_n) = \sum_J f_J x^J. 
\end{equation}
Les majorations~\eqref{eq.major.theta} des coefficients de la s\'erie~$\theta$
et l'hypoth\`ese que~$f$ est de degr\'e~$a$ entra\^{\i}nent que les
coefficients de~$F=\sum F_J x^J$ v\'erifient 
\begin{equation}\label{eq.major.L}
\begin{cases}
\ord F_J  \geq j_0 \frac{p-1}{p^2} &  \text{si $j_0a\geq j_1+\dots+j_n$,}\\
F_J = 0 & \text{sinon.}
\end{cases}
\end{equation}
Notons $L$ l'ensemble des s\'eries formelles qui v\'erifient ces
in\'egalit\'es et~$\Delta$ l'ensemble des mon\^omes $x^J$
tels que $j_0a\geq j_1+\dots+j_n$.

Soit $\psi$ l'{\og op\'erateur de Dwork\fg} sur $\Omega[[x_0,\dots,x_n]]$
donn\'e par
\begin{equation}
\psi (\sum_J u_J x^J) = \sum_J \tau(u_{pJ}) x^{J}.
\end{equation}
L'op\'erateur compos\'e~$\alpha$ de~$\psi$ et de la multiplication par~$F$
applique~$L$ dans lui-m\^eme. Soit~$A$ sa matrice (infinie)
dans la {\og base\fg}~$\Delta$.

Le point crucial de la d\'emonstration de~\textsc{Dwork}
est la nucl\'earit\'e de l'op\'erateur lin\'eaire~$\alpha^d$
qui se d\'eduit de majorations explicites pour les coefficients
de la matrice~$A$, elles-m\^emes d\'eduites de~\eqref{eq.major.L}.
En particulier, pour tout entier~$k\geq 1$, $\alpha^{dk}$ poss\`ede
une trace, d'ailleurs \'egale \`a la somme des coefficients diagonaux de sa matrice
dans la base~$\Delta$,
et l'on a
\begin{equation}
\Tr(\alpha^{k})(q^k-1)^{n+1}=q^k \Card{X(\F_{q^k})}-(q^k-1)^n.
\end{equation}
Comme l'ont remarqu\'e~\textsc{Lauder} et~\textsc{Wan}
(voir~\cite{lauder-wan2006}, th.~28), ces majorations entra\^{\i}nent
aussi que si l'on {\og tronque\fg}
l'op\'erateur~$\alpha$ en consid\'erant sa restriction~$\alpha_\nu$
au sous-espace (stable) engendr\'e par les mon\^omes~$x^J$ de~$\Delta$
tels que $j_0\leq \nu (p/(p-1))^2$, 
on obtient une congruence (pour simplifier, on ne regarde d\'esormais
que le cas $k=1$):
\begin{equation}
\Tr(\alpha_\nu)(q-1)^{n+1}\equiv q \Card{X(\F_q)}-(q-1)^n \pmod{p^\nu}
\end{equation}
Soit $A_\nu$ la matrice de l'application $\sigma$-lin\'eaire~$\alpha_\nu$; 
elle est de taille
$N=\mathrm O( (\nu a+1)^n)$. Pour calculer $\Tr(\alpha_\nu)$
modulo~$p^\nu$, il suffit de r\'eduire~$A_\nu$ modulo~$p^\nu$,
d'effectuer le produit
$A_\nu \tau(A_\nu)\dots \tau^{a-1}(A_\nu)$
puis de calculer la trace de cette matrice modulo~$p^\nu$.
Le produit de ces~$k$ matrices peut \^etre effectu\'e efficacement
en adaptant l'algorithme d'exponentiation binaire;
on n'a alors que~$\mathrm O(\log k)$ produits \`a faire.

Cela fournit un algorithme pour calculer $\Card{X(\F_q)}$
modulo~$p^\nu$. Comme $0\leq \Card{X(\F_q)}<p^{dn}$,
il suffit, pour en d\'eduire la valeur de~$\Card{X(\F_q)}$,
de choisir $\nu=nd$. La complexit\'e de l'algorithme ainsi esquiss\'e 
est~$\Otilde(p^{2n+4}a^{3n} n^{3n+5} d^{3n+7})$ en temps,
et~$\Otilde(pa^{2n}n^{2n+3}d^{2n+4})$ en espace.

\subsection{Rel\`evement canonique et moyenne arithm\'etico-g\'eom\'etrique}

Revenons au cas des courbes elliptiques.
Soit $E$ une courbe elliptique d\'efinie sur un corps fini~$\F$;
notons $q$ le cardinal de~$\F$, $p$ sa caract\'eristique;
soit $d$ l'entier tel que $q=p^d$.
Notons enfin $W$ l'anneau des vecteurs de Witt de~$\F$ et
$K$ son corps des fractions.

L'algorithme que nous d\'ecrivons ici permet, une fois encore,
de calculer~$\Card{E(\F)}$. Il est d\^u \`a~K.~\textsc{Satoh}~\cite{satoh2000},
au moins lorsque $p>3$.
La pr\'esentation que nous suivons ici tient compte
d'articles parus ult\'erieurement
dans lequel $p=2$ et~$3$ sont pris en compte
(\cite{fouquet-gaudry-harley2000},
\cite{skjernaa2003}, \cite{satoh-s-t2003} et~\cite{vercauteren-pv2001}),
ainsi que de sa variante, due \`a~\textsc{Mestre}, fond\'ee
sur la moyenne arithm\'etico-g\'eom\'etrique.

\paragraph{Rel\`evements}
Soit $\mathscr E$ un \emph{rel\`evement de~$E$},
c'est-\`a-dire une courbe elliptique sur~$W$ dont la r\'eduction
modulo~$p$ est \'egale \`a~$E$. Une telle courbe peut \^etre obtenue en
consid\'erant une \'equation de Weierstrass~\eqref{eq.weierstrass.general}
ou~\eqref{eq.weierstrass} \`a coefficients dans~$W$ 
dont la r\'eduction modulo~$p$ est une \'equation de Weierstrass de~$E$.
Un endomorphisme $\tilde u$ de~$\mathscr E$
 d\'efinit, par r\'eduction modulo~$p$, un endomorphisme de~$E$,
d'o\`u un homomorphisme d'anneaux, injectif, de~$\End(\mathscr E)$
dans~$\End(E)$.

Si $\mathscr E$ est mal choisi, $\End(\mathscr E)=\Z$.
Toutefois, si $E$ est \emph{ordinaire}, c'est-\`a-dire
que $E$ poss\`ede un point d'ordre~$p$,
un th\'eor\`eme de \textsc{Deuring}~\cite{deuring1941}
affirme qu'il existe un unique rel\`evement~$\mathscr E$
pour lequel l'homomorphisme $\End(\mathscr E)\ra\End(E)$ soit
un isomorphisme ou, cela revient au m\^eme, le seul rel\`evement
sur lequel l'endomorphisme de Frobenius~$\pi_E$ de~$E$
se rel\`eve en un endomorphisme~$\pi_{\mathscr E}$ de~$\mathscr E$.
Comme le corps~$\F$ est fini, ce rel\`evement co\"{\i}ncide
avec celui fourni par la th\'eorie, plus profonde, 
du \emph{rel\`evement canonique}, \emph{cf.}~\cite{lubin-serre-tate1964}
et, pour plus de d\'etails, \cite{messing76}.

Comme les traces de~$\pi_{\mathscr E}$ dans~$\End(\mathscr E)$
et de~$\pi_E$ dans~$\End(E)$ sont \'egales, il suffit de calculer
la premi\`ere. Celle-ci se voit au niveau des formes diff\'erentielles:
si $\omega$ est une forme diff\'erentielle invariante 
(non nulle) sur~$\mathscr E$,
$\pi_{\mathscr E}^*\omega$ est une forme diff\'erentielle invariante,
donc est proportionnelle \`a~$\omega$.
Soit $c_{\mathscr E}\in W$ 
tel que $\pi_{\mathscr E}^*\omega=c_{\mathscr E}\omega$.
Puisque l'isog\'enie duale~$\pi_{\mathscr E}^\vee$ de~$\pi_{\mathscr E}$ 
v\'erifie  la relation
\begin{equation}
\pi_{\mathscr E}^\vee\circ\pi_{\mathscr E}=q\id_{\mathscr E},
\end{equation}
$c_{\mathscr E}\neq0$ et
$(\pi_{\mathscr E}^\vee)^*\omega=qc_{\mathscr E}^{-1}\omega$.
La somme des endomorphismes $\pi_{\mathscr E}$ et~$\pi_{\mathscr E}^\vee$
est la multiplication par la trace de~$\pi_{\mathscr E}$,
c'est aussi puisque l'homomorphisme $\End(\mathscr E)\ra\End(E)$
est un isomorphisme, la multiplication par la trace
de~$\pi_E$.
Autrement dit, on a l'\'egalit\'e
\begin{equation}\label{eq.trace.c}
\Card{E(\F)}=q+1-t_E = q+1 - \left( c_{\mathscr E}+\frac q{c_{\mathscr E}}\right).
\end{equation}

Cette m\'ethode n'est pas praticable telle quelle car~$\pi_E$
est un endomorphisme de degr\'e~$q$, suppos\'e tr\`es grand. 
L'id\'ee de~\textsc{Satoh} consiste, consid\'erant
la factorisation naturelle de~$\pi_E$ en $d$~isog\'enies de degr\'e~$p$,
\`a relever ces isog\'enies et \`a calculer des constantes analogues~$c_i$
dont~$c_{\mathscr E}$ sera le produit.
En fait, une telle factorisation est un point essentiel
de la construction m\^eme de la courbe~$\mathscr E$,
suivant la m\'ethode d\'ecrite dans~\cite{lubin-serre-tate1964}.

Soit $\sigma$ l'automorphisme de~$\F$ donn\'e par $x\mapsto x^p$;
comme $q=p^d$, on a $\sigma^d=\id$.
Si $E$ est une courbe elliptique sur~$\F$,
on notera $E^{(\sigma)}$ 
la courbe obtenue en appliquant l'automorphisme~$\sigma$
\`a ses coefficients; elle est li\'ee \`a~$E$ par une isog\'enie 
$\phi\colon E\ra E^\sigma$ de degr\'e~$p$ qui correspond \`a l'\'el\'evation
des coordonn\'ees des points \`a la puissance~$p$.  
On a ainsi une suite d'isog\'enies (de degr\'e~$p$)
\begin{equation}
  E \xrightarrow{\phi}  E^{(\sigma)} \xrightarrow{\phi} E^{(\sigma^2)}
\xrightarrow\phi\dots\xrightarrow{\phi} E^{(\sigma^d)}=E .
\end{equation}

Notons encore $\sigma$ l'unique automorphisme de~$W$ dont la r\'eduction
modulo~$p$ est l'automorphisme~$\sigma$ de~$\F$; on a aussi
$\sigma^d=\id$.
\`A la suite de~\cite{lubin-serre-tate1964},
\textsc{Satoh} propose de chercher~$\mathscr E$ de sorte que l'on
ait encore une suite d'isog\'enies de degr\'e~$p$:
\begin{equation}
 \mathscr E \xrightarrow{\phi}  \mathscr E^{(\sigma)} \xrightarrow{\phi^{(\sigma)}} \mathscr E^{(\sigma^2)}
\xrightarrow{\phi^{(\sigma^2)}}\dots\xrightarrow{\phi^{(\sigma^{d-1})}} \mathscr E^{(\sigma^d)}=\mathscr E .
\end{equation}
En effet, \'etant donn\'ee une telle courbe~$\mathscr E$, l'endomorphisme
compos\'e $\phi^{(\sigma^{d-1})}\circ\dots\circ\phi^{(\sigma)}\circ\phi$
rel\`eve l'endomorphisme~$\pi_E$ de~$E$, si bien que $\mathscr E$
est le rel\`evement canonique de~$E$.
Autrement dit, on va chercher un rel\`evement~$\mathscr E$
de~$E$ tel que le couple $(\mathscr E, \mathscr E^{(\sigma)})$
soit un point de la courbe modulaire~$X_0(p)$.

\paragraph{R\'esolution de l'\'equation modulaire}
Le polyn\^ome modulaire~$\Phi_p$,
\'equation de la courbe~$X_0(p)$ dans~$\P^1\times\P^1$,
v\'erifie une congruence modulo~$p$, due \`a \textsc{Kronecker}:
\begin{equation}
\Phi_p (X,Y) \equiv (X^p-Y)(Y^p-X) \pmod p.
\end{equation}
Lorsque l'invariant $j_E$ de la courbe~$E$
n'appartient pas au corps~$\F_{p^2}$,  
l'\'equation {\og sesquipolynomiale\fg}
$\Phi_p(j,\sigma(j))$, o\`u l'inconnue~$j$
appartient \`a l'ensemble des \'el\'ements de~$W$ dont la r\'eduction
est l'invariant~$j_E$,
est alors justiciable d'une variante de la m\'ethode de Newton,
d'o\`u une construction du rel\`evement canonique de la courbe~$E$
dans ce cas.
\`A chaque it\'eration, 
l'\'equation lin\'earis\'ee
est de la forme, dite \emph{\'equation d'Artin-Schreier}:
\begin{equation}
a\sigma(x)+bx=c.
\end{equation}
Toutefois, le calcul effectif de~$\sigma$ est d\'elicat.
Les algorithmes efficaces de~\textsc{Lercier}, \textsc{Lubicz}
puis~\textsc{Harley} 
pour les r\'esoudre leur ont permis d'impl\'ementer effectivement
cette approche
(voir~\cite{lercier-lubicz2003,harley2002}
et le chapitre~12 de~\cite{cohen2006}).
Quant \`a \textsc{Satoh}, il proc\'eda un peu diff\'eremment
en prenant pour inconnue, non pas seulement la courbe~$\mathscr E$, mais 
toutes les courbes~$\mathscr E_i=\mathscr E^{(\sigma^i)}$ simultan\'ement.
Les $d$ invariants $j_i=j(\mathscr E_i)$, pour $0\leq i\leq d-1$,
sont li\'es par le syst\`eme
d'\'equations polynomiales 
\begin{equation}
\Phi_p(j_0,j_1) =\Phi_p(j_1,j_2)=\dots=\Phi_p(j_{d-2},j_{d-1})=\Phi_p(j_{d-1},j_0)
\end{equation}
dont une r\'esolution par la m\'ethode de Newton classique est possible.

\paragraph{It\'eration de quotients}
Avec les notations pr\'ec\'edentes, les courbes~$\mathscr E_0$
et~$\mathscr E_1$ sont li\'ees par une isog\'enie 
$\phi\colon\mathscr E_0\ra\mathscr E_1$ dont la r\'eduction modulo~$p$,
le frobenius relatif, est une isog\'enie ins\'eparable.
Soit $\phi\colon\mathscr E_0\ra\mathscr E_1$ un rel\`evement arbitraire 
de~$\phi\colon E\ra E^{(\sigma)}$.
Le sous-sch\'ema de~$p$-torsion, $\mathscr E_0[p]$, poss\`ede
un sous-groupe de rang~$p$ connexe canonique, $G_0$, le quotient
$\mathscr E_0[p]/G_0$ \'etant un $W$-sch\'ema en groupes \'etale.
(Localement pour la topologie \'etale, $G_0$ est isomorphe 
au sch\'ema en groupes~$\mu_p$.)
De plus, $\mathscr E_1$ est isomorphe \`a~$\mathscr E_0/G_0$
et $\phi$ s'identifie \`a l'isog\'enie $\mathscr E_0\ra \mathscr E_0/G_0$.
Ainsi, $\mathscr E_1$ est d\'etermin\'ee par~$\mathscr E_0$.

Cette analyse faite,
choisissons une courbe elliptique arbitraire~$\mathscr E_0$
relevant~$E$, soit $G_0$ le plus grand sous-sch\'ema en groupes
connexe de~$\mathscr E_0[p]$,
soit $\mathscr E_1$ le quotient de~$\mathscr E_0$
par~$G_0$
et soit $\phi_0\colon\mathscr E_0\ra\mathscr E_1$ l'isog\'enie
canonique. La courbe elliptique~$\mathscr E_1$ rel\`eve
la courbe ordinaire~$E^{(\phi)}$. It\'erons le proc\'ed\'e.
On en d\'eduit une suite~$(\mathscr E_n)$ de courbes elliptiques
li\'ees par des isog\'enies de degr\'e~$p$,
$\phi_i\colon\mathscr E_n\ra\mathscr E_{n+1}$,
dont le noyau, $G_n$, est le plus grand sous-sch\'ema en groupes
connexe de~$\mathscr E_n[p]$.

Comme l'application $x\mapsto x^p$ est $p$-adiquement contractante,
les sous-suites~$(\mathscr E_{dn+j})$ sont convergentes
dans la courbe modulaire. En particulier, les courbes~$\mathscr E_{dn}$
convergent vers le rel\`evement canonique de~$E$
lorsque $n\ra\infty$.

Une fa\c{c}on de le d\'emontrer consisterait \`a utiliser
le sch\'ema formel des rel\`evements de la courbe~$E$
(\cite{lubin-serre-tate1964,messing76}), 
lequel s'identifie au sch\'ema formel des rel\`evements de son groupe~$p$-divisible.
La th\'eorie des coordonn\'ees canoniques (qui s'appellent~$q$
d'ordinaire, mais que je noterai~$z$ ici) montre 
en effet qu'il est isomorphe au groupe
multiplicatif formel sur~$W$; l'\'el\'ement neutre correspond au rel\`evement
canonique de la courbe~$E$. 
Par ailleurs, les isog\'enies que nous consid\'erons (ou plut\^ot
les compos\'ees de $d$ isog\'enies successives) ont pour
effet d'\'elever \`a la puissance~$q$ la coordonn\'ee canonique~$z$
essentiellement parce que l'\'el\'evation \`a la puissance~$q$ 
dans le groupe multiplicatif formel a pour noyau le sch\'ema
en groupes~$\mu_q$. Comme $\abs{z-1}<1$, l'it\'eration
de l'application $z\mapsto z^q$ fournit une suite qui converge vers~$1$.
On observe que la convergence n'est que lin\'eaire.

Cette construction se g\'en\'eralise ainsi \emph{mutatis mutandis}
au cas des vari\'et\'es ab\'eliennes
ordinaires (voir aussi le chapitre~2 de~\cite{carls2004}).
Citons aussi l'article~\cite{kohel2003} qui propose
des variantes de ces algorithmes o\`u la courbe
modulaire correspondant \`a l'invariant~$j$
est remplac\'ee par la courbe~$X_0(N)$, surtout
lorsque celle-ci est de genre~$0$.

Dans le cas des courbes elliptiques,
son int\'er\^et fut mis en \'evidence dans l'article~\cite{vercauteren-pv2001}
de \textsc{Vercauteren}, \textsc{Preenel}
et~\textsc{Vandewalle}, 
parce que sa consid\'eration am\'eliore la complexit\'e en espace de l'algorithme
de~\textsc{Satoh}.
Toutefois, dans cet article, l'invariant~$j_{n+1}$ est calcul\'e en appliquant
une it\'eration de l'algorithme de Newton \`a l'\'equation
modulaire~$\Phi_p(j,j_n)=0$, l'inconnue \'etant~$j$,
le tout \'etant \'ecrit \`a l'envers, car ces auteurs raisonnent,
ainsi que~\textsc{Satoh}, sur l'isog\'enie de Verschiebung.

\paragraph{Fin de l'algorithme}
Pour calculer explicitement sous forme
d'une \'equation de Weierstrass~\eqref{eq.weierstrass.general}
ces courbes~$\mathscr E_n$
et les isog\'enies~$\phi_n$, on peut utiliser
les formules de \textsc{V\'elu}~\cite{velu1971},
voir aussi~\cite{carls2004}, th.~3.3.1.

Supposons donc la courbe~$\mathscr E$
et l'isog\'enie~$\phi$ explicit\'ees et montrons
comment terminer le calcul du cardinal de~$E(\F)$.
Soit $\omega$ une base du $W$-module des formes diff\'erentielles invariantes
sur~$\mathscr E$. Il existe  un unique \'el\'ement $c\in W$ 
tel que $\phi^* \omega^{(\sigma)}=c\omega$, et les formules
de~\textsc{V\'elu} permettent d'ailleurs de le calculer ais\'ement.
Pour $i\in\{0,\dots,d-1\}$, on a donc
$(\phi^{(\sigma^i)})^*\omega^{(\sigma^{i+1})}=\sigma^i(c)\omega^{(\sigma^i)}$,
si bien que
\begin{equation}
 c_{\mathscr E} = c \sigma(c) \dots \sigma^{d-1}(c)=N_{K/\Q_p}(c)
\end{equation}
est la norme de~$c$ dans l'extension~$\Q_p\subset K$.

Si l'on ne conna\^{\i}t qu'une valeur approch\'ee de~$c$, disons
modulo~$p^N$, on en d\'eduit une valeur de~$c_{\mathscr E}$ modulo~$p^N$,
donc, via la formule~\eqref{eq.trace.c}, la valeur de~$t_{E}$
modulo~$p^{N-d}$, la perte de pr\'ecision \'etant due au fait
que $c\equiv 0\pmod p$.
C'est pour pallier cet inconv\'enient
que \textsc{Satoh} utilise
l'isog\'enie duale~$\psi$, dite \emph{Verschiebung}; comme $E$
est ordinaire, $\psi$ est  s\'eparable et 
la constante analogue~$c$  qui intervient dans le calcul de la trace
est une unit\'e de~$W$.

Une fois calcul\'e $t_E$ modulo~$p^N$, l'in\'egalit\'e de Hasse
$\abs{t_E}\leq 2\sqrt q$
permet d'en d\'eduire~$t_E$ si $p^N>4\sqrt q$, c'est-\`a-dire si
$N>\frac d2$, voire $\frac{d}2+2$ si $p=2$.

D\'ecrivons rapidement la complexit\'e de l'algorithme obtenu.
Comme il faut calculer~$c$ modulo~$p^{N}$, o\`u $N=\mathrm O(d)$,
et que la pr\'ecision augmente de~$1$ \`a chaque it\'eration,
le nombre d'it\'erations \`a effectuer
est~$\mathrm O(d)$. Par ailleurs, la taille des objets
\`a manipuler, des \'el\'ements de $(\Z/p^N\Z)[x]/(f(x))$
dont la multiplication n\'ecessite~$\Otilde(Nd\log p)$ op\'erations,
d'o\`u une complexit\'e de~$\Otilde (d^2\log p)$ op\'erations
pour calculer~$c$. Comme l'a remarqu\'e~\textsc{Harley}~\cite{harley2002},
les techniques efficaces pour calculer les r\'esultants permettent
de calculer la norme de~$c$ modulo~$p^N$ en~$\Otilde(Nd\log p)$
op\'erations. Par suite, la complexit\'e de l'algorithme 
de calcul de~$\Card{E(\F)}$ suivant la m\'ethode de~\textsc{Satoh}
est $\Otilde(d^2\log p)$. 
Suite \`a~\cite{vercauteren-pv2001},
la complexit\'e en espace est du m\^eme ordre.

Il convient de remarquer qu'elle n'est polynomiale en~$\log q$
que si la caract\'eristique~$p$ du corps est fixe.
Toutefois, pour les applications cryptographiques, on choisit
souvent $p=2$ et l'algorithme ainsi obtenu est tr\`es efficace.
En 2002, \textsc{Fouquet}, \textsc{Gaudry} et~\textsc{Harley}
ont ainsi pu calculer 
le nombre de points d'une courbe
elliptique sur un corps fini de cardinal~${2^{8\,009}}$
en environ~$300$~h de calcul;
le stockage des donn\'ees a n\'ecessit\'e plus de~15~Go 
(voir~\cite{fouquet-gaudry-harley2000})!

\paragraph{Moyenne arithm\'etico-g\'eom\'etrique}
La moyenne arithm\'etico-g\'eom\'etrique $\mathrm M(a,b)$
de deux nombres r\'eels, disons
strictement positifs, $a$ et~$b$, est la limite commune
des deux suites~$(a_n)$ et~$(b_n)$ d\'efinies par les relations
de r\'ecurrence
\begin{equation}\label{eq.agm.0}
a_{n+1}=\frac{a_n+b_n}2, \quad b_{n+1}=\sqrt{a_nb_n},
\end{equation}
et l'initialisation $a_0=a$, $b_0=b$.
Ainsi que l'a d\'ecouvert C.-F.~\textsc{Gauss} en 1799,
elle est reli\'ee aux int\'egrales elliptiques par les formules
\begin{equation}\label{eq.gauss.agm}
 \frac{\pi}{\mathrm M(a,b)}
= 2\int_{0}^{\pi/2} \frac{dt}{\sqrt{a^2\cos^2t+b^2\sin^2t}}
= \int_0^\infty \frac{dx}{\sqrt{x(x+a^2)(x+b^2)}}
\end{equation}
(changement de variables $x=b^2\tan^2 t$)
et permettent donc de calculer tr\`es rapidement
les p\'eriodes de la courbe elliptique d'\'equation
\begin{equation}
y^2=x(x+a^2)(x+b^2).
\end{equation}
(L'autre p\'eriode est obtenue en consid\'erant la courbe
tordue, d'\'equation $y^2=x(x-a^2)(x-b^2)$, et
fait intervenir $\mathrm M(a+b,a-b)$.)
La raison d'\^etre de la formule~\eqref{eq.gauss.agm}
est l'existence d'une isog\'enie de degr\'e~$2$ entre la courbe elliptique
pr\'ec\'edente et celle d'\'equation
\begin{equation}
y^2=x\big( x+ \big(\frac{a+b}2\big)^2 \big) \big( x+ab\big).
\end{equation}
Cette isog\'enie fournit un changement de variables dans l'int\'egrale elliptique
qui ram\`ene \`a une int\'egrale similaire o\`u $a$ et~$b$ sont remplac\'es
par~$a_1$ et~$b_1$.
Les suites~$(a_n)$ et~$(b_n)$ apparaissent ainsi comme un proc\'ed\'e
simple pour calculer une suite d'isog\'enies de degr\'e~$2$ particuli\`eres,
en l'occurrence la seule fournissant une suite de courbes
elliptiques dont les points d'ordre~$2$ sont tous d\'efinis sur~$\R$.

En~2000, J.-F.~\textsc{Mestre}, voir~\cite{mestre2000},
avait montr\'e comment une variante $2$-adique de la moyenne
arithm\'etico-g\'eom\'etrique permet de calculer
efficacement rel\`evement canonique et nombre de points
lorsque $p=2$.   
L'algorithme~\textsc{agm} qui en r\'esulte est essentiellement
\'equivalent (mais ant\'erieur !) \`a la modification
par~\textsc{Vercauteren}~\cite{vercauteren-pv2001}
de celui de~\textsc{Satoh};
son impl\'ementation est cependant tr\`es ais\'ee 
et son ex\'ecution plus rapide.
\textsc{Gaudry}, \textsc{Harley}, \textsc{Lercier},
\textsc{Lubicz} le mirent en \oe uvre
fin~2002 dans des corps de tr\`es grand degr\'e,
le record semblant \^etre d\'etenu par \textsc{Harley} 
qui a pu calculer un tel cardinal~$\Card{E(\F_q)}$,
pour $q=2^{130\,020}$. (Un tel corps poss\`ede une base normale engendr\'ee
par une somme de Gauss, ce qui permet d'acc\'el\'erer tr\`es notablement 
certains algorithmes.)

Profitons-en pour donner des formules.

Dans la suite, nous raisonnerons en fait sur la quantit\'e $\xi=a/b$,
dont le carr\'e est l'invariant de Legendre de la courbe elliptique en question.

Dans ce paragraphe, supposons que $p=2$ et que $q=p^d$, pour $d\geq 1$.
Soit $E$ une courbe elliptique ordinaire sur~$\F_q$, donn\'ee par une \'equation
\begin{equation}\label{eq.weierstrass.2}
y^2+xy=x^3+c
\end{equation}
o\`u $c\in\F^*$ est l'inverse de l'invariant~$j$. (Ces courbes
sont celles dont la $4$-torsion est d\'efinie sur~$\F_q$;
on s'y ram\`ene par une torsion quadratique.)
Comme dans le cas complexe, le calcul de la moyenne
arithm\'etico-g\'eom\'etrique n\'ecessite un choix de racines carr\'ees ;
pour tout $t\in W$ tel que $t\equiv 1\pmod 8$, on note $\sqrt t$
l'unique \'el\'ement de~$W$ congru \`a~$1$ modulo~$4$.

Soit $\xi_4$ un \'el\'ement de~$W$ tel que $\xi_4\equiv 1+8c\pmod {16}$.
On d\'efinit alors une suite $(\xi_n)$ par r\'ecurrence 
en posant, pour $n\geq 4$,
\begin{equation}
\xi_{n+1} = \frac{1+\xi_n}{2\sqrt{\xi_n}}.
\end{equation}
Pour $0\leq i<d$, la  suite~$(\xi_{dn+i})$ converge
vers un \'el\'ement~$\xi^*_i$ dont le carr\'e est l'invariant
de Legendre de la courbe~$\mathscr E_i=\mathscr E^{(\sigma ^i)}$,
$\mathscr E$ d\'esignant le  rel\`evement canonique de~$E$.
En outre, si l'on pose
\begin{equation}
\mu_n=\frac{2\xi_n}{1+\xi_n}
\quad\text{et}\quad
 t_n = \mu_n \sigma(\mu_n) \dots \sigma^{d-1}(\mu_n)= N_{K/\Q_2}(\mu_n),
\end{equation}
on a  pour tout entier~$n\geq 0$
la congruence
\begin{equation}
q+1-\Card{E(\F)}=\Tr(\pi_E) \equiv t_{n}+\frac q{t_n} \pmod{2^{n}},
\end{equation}
d'o\`u le calcul de~$\Card{E(\F)}$ d\`es que $2^n>4\sqrt q$,
c'est-\`a-dire $n>\frac 12d+2$.

%

\paragraph{G\'en\'eralisation aux courbes de genre sup\'erieur}
En 1836, \textsc{Richelot}  a d\'emontr\'e l'existence d'un algorithme
permettant de calculer
les {\og int\'egrales hyperelliptiques de genre~$2$\fg},
g\'en\'eralisant ainsi la formule~\eqref{eq.gauss.agm}.
Comme dans le cas elliptique, ce th\'eor\`eme s'interpr\`ete 
en termes d'isog\'enies de surfaces ab\'eliennes,
voir~\cite{bost-mestre1988}. Dans sa lettre~\cite{mestre2000},
\textsc{Mestre} proposait d'utiliser ces relations de r\'ecurrence
dans le cas~$2$-adique pour calculer le nombre de points
d'une courbe de genre~$2$ d\'efinie sur un corps fini de caract\'eristique~$2$.

Dans~\cite{mestre2002}, il reprend ce sujet en genre sup\'erieur
\`a l'aide des formules de duplication des fonctions~$\theta$.
Elles permettent de calculer par un proc\'ed\'e it\'eratif 
(convergeant ver la \emph{moyenne de Borchardt})
un d\'eterminant~$2\times 2$
form\'e \`a l'aide des p\'eriodes r\'eelles d'une courbe hyperelliptique
de genre~$2$ d\'efinie sur~$\R$
dont les points de Weierstrass sont r\'eels;
l'initialisation de l'algorithme utilise les formules de \textsc{Thomae}
(voir par exemple~\cite{mumford1984}, chap.~3A, \S8).

\textsc{Mestre} montre que cet algorithme est d'une grande
utilit\'e pour calculer le nombre de points d'une courbe hyperelliptique
\emph{ordinaire} sur un corps fini.
Pr\'ecis\'ement, \'etant donn\'ee une courbe~$C$ de genre~$2$ sur un corps
fini~$\F_q$, de caract\'eristique~$2$,
cet algorithme, convenablement interpr\'et\'e dans l'anneau~$W$
des vecteurs de Witt,
permet de calculer le rel\`evement~$\mathscr C$ de~$C$
sur~$W$  dont la jacobienne~$J_{\mathscr C}$ est le rel\`evement
canonique de celle, $J_C$, de~$C$.
En outre, il fournit le d\'eterminant de l'endomorphisme
relevant celui de Verschiebung~$\psi_q$ agissant sur le $W$-module
des formes diff\'erentielles globales de~$\mathscr C$.
Autrement dit, des quatre valeurs propres du frobenius
$\pi_1,\dots,\pi_4$
on peut d\'eterminer le produit~$\pi_1\pi_2$ des deux
qui sont inversibles modulo~$2$.

Cette derni\`ere partie de l'algorithme, \`a savoir l'obtention
du produit~$\alpha=\pi_1\dots\pi_g$
des valeurs propres $\pi_1,\dots,\pi_{2g}$
qui sont des unit\'es~$2$-adiques,
se g\'en\'eralise aux courbes hyperelliptiques ordinaires de tout genre~$g$.
Elle a en outre \'et\'e \'etendue aux courbes de genre~$3$
non hyperelliptiques par~\textsc{Ritzenthaler}~\cite{ritzenthaler2004};
au lieu des points de Weierstrass qui apparaissent dans les formules
de \textsc{Thomae}, il fait usage des bitangentes.
Par contre, l'algorithme ne calcule pas le rel\`evement
canonique de la jacobienne, pour la bonne raison que ce n'est
pas forc\'ement une jacobienne si $g\geq 4$.\footnote{Pour
\^etre pr\'ecis, ce r\'esultat n'est prouv\'e dans~\cite{dwork-ogus1986}
que sous l'hypoth\`ese que la caract\'eristique du corps est
diff\'erente de~$2$...}

Une fois obtenu~$\alpha$, il n'est pas toujours
possible d'en d\'eduire les~$\pi_i$, au moins lorsque $g\geq 4$.
\textsc{Mestre} donne par exemple
un exemple de vari\'et\'es ab\'eliennes de dimension~$4$ sur~$\F_2$,
ordinaires,
non g\'eom\'etriquement isog\`enes et dont les invariants~$\alpha$
sont \'egaux.
Posons toutefois $\beta=\alpha+q^g\alpha^{-1}$
et consid\'erons, \`a la suite de~\textsc{Mestre}, le polyn\^ome
minimal~$P$ de~$\beta$ sur~$\Z$ (calcul\'e par exemple
\`a l'aide d'une version~$2$-adique de l'algorithme~LLL).
Si $P$ est de degr\'e~$2^{g-1}$, 
ce qui arrive si par exemple
la jacobienne~$J_C$ de~$C$ est simple, les racines de~$P$
sont les \'el\'ements 
\[ 
q^{\Card{\complement I}} \frac{\prod_{i\in I}\pi_i}{\prod_{i\not\in I}\pi_i }
 + q^{\Card{I}} \frac{\prod_{i\not\in I}\pi_i}{\prod_{i\in I}\pi_i} .\]
Cela permet de d\'eterminer les~$\pi_i$ au signe pr\`es,
d'o\`u le cardinal de~$J_C(\F_q)$ \`a une ambigu\"{\i}t\'e finie pr\`es.
 \`A l'oppos\'e, si $P$ est degr\'e~$1$, c'est-\`a-dire si $\beta\in\Z$,
alors $J_C$ est g\'eom\'etriquement isog\`ene \`a une puissance d'une courbe
elliptique. Ces deux remarques
permettent de d\'em\^eler  la situation lorsque $g=2$ ou~$3$.

En genre~$2$, l'algorithme a \'et\'e \'etudi\'e et impl\'ement\'e
par~\textsc{Lercier} et~\textsc{Lubicz}~\cite{lercier-lubicz2006}; sa complexit\'e est~$\Otilde(n^2)$, tant en espace
qu'en temps. Il leur a permis de calculer en quelques jours
le cardinal de courbes de genre~$2$ sur un corps \`a $2^{32\,770}$ \'el\'ements,
et de genre~$3$ sur un corps \`a $2^{4\,098}$~\'el\'ements.
En genre~$3$,
cette m\'ethode a permis \`a 
\textsc{Ritzenthaler} (voir~\cite{ritzenthaler2004})
de calculer en deux semaines
le nombre de points  rationnels
d'une courbe
quartique plane (une courbe non hyperelliptique de genre~$3$)
sur un corps de cardinal~$2^{5\,002}$.

Tr\`es r\'ecemment, \textsc{Carls} et~\textsc{Lubicz}
ont annonc\'e dans~\cite{carls-lubicz2007}
l'existence d'un algorithme (quasi-quadratique 
en l'exposant de la caract\'eristique du corps) permettant
de calculer le nombre de points d'une courbe hyperelliptique
ordinaire, g\'en\'eralisant ainsi
en toute caract\'eristique les algorithmes
de~\textsc{Satoh} et~\textsc{Mestre}.
Cet algorithme est fond\'e sur l'existence, 
due \`a~\textsc{Carls}~\cite{carls2007},
de structures th\^eta
canoniques sur le rel\`evement canonique d'une vari\'et\'e ab\'elienne ordinaire.

\subsection{Cohomologie de Monsky--Washnitzer}

En 2001, K.~\textsc{Kedlaya}~\cite{kedlaya2001} a propos\'e  un algorithme
calculant le nombre de points de courbes hyperelliptiques.
Le r\^ole principal n'est plus tenu par la jacobienne de la courbe,
comme dans  les algorithmes pr\'ec\'edents, mais par sa cohomologie
de Monsky--Washnitzer dont \textsc{Kedlaya} montre qu'elle
peut \^etre calcul\'ee efficacement, au moins si la caract\'eristique
du corps est petite.

Soit $\F$ un corps fini, de caract\'eristique~$p$, de cardinal~$q$.
Soit $W$ l'anneau des vecteurs de Witt de~$\F$ et soit $K$
son corps des fractions. 
Notons $\sigma$ l'automorphisme de Frobenius sur~$\F$, son
rel\`evement \`a~$W$ et son extension \`a~$K$; si $q=p^d$, on a $\sigma^d=\id$.

\paragraph{D\'efinition de la cohomologie de Monsky--Washnitzer}
Rappelons ce dont il s'agit, en renvoyant \`a~\cite{monsky-w68,berthelot1996,
berthelot1988,vanderput1988,fresnel-vanderput2004,lestum2007}
pour plus de d\'etails.
Consid\'erons une vari\'et\'e alg\'ebrique affine et lisse~$X$ d\'efinie 
sur~$\F$.
Le but est de d\'efinir une sorte de cohomologie de De Rham de~$X$
munie d'un endomorphisme de Frobenius induit par celui de~$X$,
de sorte \`a ce que la fonction z\^eta de~$X$ se calcule par une formule
de Lefschetz. Le premier pas consiste \`a relever~$X$ en caract\'eristique~$0$.

D'apr\`es un th\'eor\`eme d'\textsc{Elkik}~\cite{elkik1973},
il existe un $W$-sch\'ema affine et lisse~$\mathscr X$
dont la r\'eduction modulo~$p$ est \'egale \`a~$X$.
Soit $\mathscr X\hra\A^n_W$ une immersion ferm\'ee de~$\mathscr X$
dans l'espace affine sur~$W$, et soit $f_1,\dots,f_m\in W[x_1,\dots,x_n]$
des g\'en\'erateurs de l'id\'eal de~$\mathscr X$.
Si $A$ d\'esigne l'alg\`ebre $W[x_1,\dots,x_n]/(f_1,\dots,f_m)$ de~$\mathscr X$,
on a donc un isomorphisme entre~$A/pA$ et l'anneau~$\bar A$ 
des fonctions de~$X$.

Notons $\Omega^*_{A/W}$ le complexe de De Rham de~$A$;
le $A$-module $\Omega^1_{A/W}$ est engendr\'e par 
des \'el\'ements~$\dx_1,\dots,\dx_n$,
li\'es par les relations 
$\mathrm df_j=\sum_{i=1}^n \frac{\partial f_j}{\partial x_i} \dx_i=0$;
enfin, $\Omega^k_{A/W}=\bigwedge^k \Omega^1_{A/W}$.
La cohomologie de ce complexe est assez pathologique, m\^eme si $\mathscr X$
est la droite affine ;  en revanche, apr\`es tensorisation par~$K$, 
on obtient 
d'apr\`es un th\'eor\`eme de~\textsc{Grothendieck}~\cite{grothendieck1966}
la cohomologie de De Rham de~$\mathscr X_K$.

Toutefois, l'endomorphisme (semi-lin\'eaire) de Frobenius de~$\bar A$,
$\phi\mapsto x\mapsto x^p$, n'a en g\'en\'eral pas de rel\`evement \`a~$A$;
dans le cas des courbes elliptiques ordinaires, il aurait fallu
avoir choisi pr\'ecis\'ement le rel\`evement canonique.

Introduisons alors la $W$-alg\`ebre $W[x_1,\dots,x_n]^\dagger$ 
des s\'eries formelles~$f=\sum a_mx^m\in W[[x_1,\dots,x_n]]$
qui convergent dans un polydisque de rayon~$>1$, autrement
dit, dont les coefficients $a_m$ v\'erifient une in\'egalit\'e
de la forme $\abs{a_m}\leq C \rho^{\abs m}$, avec $\rho<1$.
Par d\'efinition, l'espace-dague~$A^\dagger$, 
dit encore \emph{compl\'et\'e faible},
de l'alg\`ebre~$A$
est le quotient de l'alg\`ebre~$W[x_1,\dots,x_n]^\dagger$
par l'id\'eal engendr\'e par les polyn\^omes~$f_i$.
C'est un rel\`evement de l'anneau de~$X$, au sens o\`u 
$A^\dagger/pA^\dagger\simeq \bar A$.
\`A isomorphisme pr\`es, il ne d\'epend que de la r\'eduction modulo~$p$
de l'alg\`ebre~$A$, c'est-\`a-dire que de~$X$.
L'introduction de~$A^\dagger$ permet en outre l'existence d'un rel\`evement
semi-lin\'eaire, not\'e~$\phi$,
de l'endomorphisme de Frobenius de~$\bar A$ \`a~$A^\dagger$.

Le complexe de De Rham de l'alg\`ebre~$A^\dagger$
est le complexe
$\Omega^*_{A^\dagger/W}$ des formes diff\'erentielles \emph{surconvergentes},
donn\'e par
\[ \Omega^1_{A^\dagger/W}= \Omega^1_{A/W}\otimes_A A^\dagger,
\qquad \Omega^k_{A^\dagger /W}=\bigwedge^k \Omega^1_{A^\dagger/W}, \]
avec la diff\'erentielle \'evidente; c'est un rel\`evement du complexe
de De Rham de~$X$. 
La cohomologie de Monsky--Washnitzer de~$X$
est alors d\'efinie par
\begin{equation}
H^i_\MW(X/K) = H^i( \Omega^*_{A^\dagger/W})\otimes_W K. 
\end{equation}
L'action de l'endomorphisme~$\phi$ 
de~$A^\dagger$ induit 
un endomorphisme semi-lin\'eaire bijectif de~$H^i_\MW(X/K)$,
c'est-\`a-dire que l'on a
\begin{equation}
\phi( a \omega ) = \sigma(a) \phi(\omega), \qquad\text{pour $a\in K$, $\omega\in H^i_{\MW}(X/K)$}, 
\end{equation}

M\^eme si ce n'est pas \'evident sur leur d\'efinition,
ces $K$-espaces vectoriels, de m\^eme que l'action de~$\phi$,
ne d\'ependent que de~$X$,
et sont fonctoriels en~$X$
(la tensorisation par~$K$ est n\'ecessaire \`a ce point).

Ils sont nuls pour $i<0$ ou~$i>\dim X$, 
de dimension finie (r\'esultat d\^u ind\'ependamment 
\`a \textsc{Berthelot}~\cite{berthelot1997},
\textsc{Kedlaya}~\cite{kedlaya2006}
et
\textsc{Mebkhout}~\cite{mebkhout1997})
et donnent lieu \`a une formule des traces de Lefschetz:
\begin{equation}\label{eq.rig.lefschetz}
\Card{X(\F)} = \sum_{i=0}^{\dim X} (-1)^i q^{\dim X} \Tr(q^{\dim X}F_X^{-1}|H^i_{\MW}(X/K)),
\end{equation}
o\`u $F_X=\phi^d$ est aussi induit par l'endomorphisme de Frobenius de~$X$,
donn\'e par l'\'el\'evation \`a la puissance~$q$ sur le faisceau~$\mathscr O_X$.

La d\'efinition a \'et\'e \'etendue par P.~\textsc{Berthelot}
au cas des vari\'et\'es alg\'ebriques g\'en\'erales sur~$\F$,
mais nous ne consid\'ererons ici que celui de vari\'et\'es
lisses.
Les espaces de cohomologie rigide $H^i_\rig(X/K)$
qu'il d\'efinit fonctoriellement sont des $K$-espaces
vectoriels, munis d'un endomorphisme semi-lin\'eaire
bijectif {\og de Frobenius\fg}~$\phi$ (un $\sigma$-isocristal). 
Ils sont nuls pour $i<0$ ou~$i>2\dim X$, de dimension finie, et
donnent lieu \`a une formule des traces
analogue \`a la pr\'ec\'edente, si ce n'est que la somme va de~$i=0$
\`a~$i=2\dim X$. 

L'apparition de~$q^{\dim X}F_X^{-1}$ dans la formule~\eqref{eq.rig.lefschetz}
rappelle qu'il s'agit de cohomologie sans support;
\textsc{Berthelot} a aussi d\'efini une cohomologie rigide \`a supports 
compacts, reli\'ee dans le cas lisse \`a la cohomologie rigide par une dualit\'e
de Poincar\'e~\cite{berthelot1997b}.

Lorsque $X$ est affine et lisse, ses espaces de cohomologie
rigide co\"{\i}ncident avec ceux d\'efinis par~\textsc{Monsky} et~\textsc{Washnitzer}.
Par ailleurs, lorsque $X$ est propre et lisse, sa cohomologie cristalline
$H^i_\cris(X/W)$ d\'efinit un $W$-r\'eseau de~$H^i_\rig(X/K)$,
stable par~$\phi$. 
Lorsque, de plus, $X$ est la r\'eduction
modulo~$p$ d'un $W$-sch\'ema propre et lisse~$\mathscr X$,
$H^i_\cris(X/W)$ s'identifie \`a la cohomologie de De Rham de~$\mathscr X$,
d\'efinie comme l'hypercohomologie du complexe 
de De Rham~$\Omega^*_{\mathscr X}$.  

Plus g\'en\'eralement, soit $\mathscr X^*$ un $W$-sch\'ema propre et lisse,
soit $\mathscr Y$ un diviseur de Cartier relatif de~$\mathscr X^*$
lisse sur~$W$ (voire \`a croisements normaux stricts)
et supposons que $\mathscr X=\mathscr X^*\setminus\mathscr Y$.
Notons $X^*$, $Y$, $X$ les r\'eductions de~$\mathscr X^*$, $\mathscr Y$
et~$\mathscr X$.
Par un th\'eor\`eme de~\textsc{Baldassarri} 
et~\textsc{Chiarellotto}~\cite{baldassarri-chiarellotto1994},
la cohomologie de De Rham alg\'ebrique $H^i_{\dR}(\mathscr X_K/K)$
et la cohomologie rigide~$H^i_\rig(X/K)$ sont des~$K$-espaces
vectoriels isomorphes
(voir aussi~\cite{chiarellotto-ls2002} o\`u
il est montr\'e que l'isomorphisme construit pr\'ec\'edemment
est compatible aux poids).
Si, de plus, $\mathscr X$ est affine,
d'anneau~$A$, sa cohomologie rigide $H^i_\rig(X/K)=H^i_\MW(X/K)$, d\'efinie
comme la cohomologie du complexe $\Omega^*_{A^\dagger/W}\otimes K$,
est donc \'egale \`a la cohomologie du complexe $\Omega^*_{A/W}\otimes K$.
En outre, d'apr\`es un th\'eor\`eme de~\textsc{Atiyah} et~\textsc{Hodge}
(\cite{atiyah-hodge1955}, voir aussi~\cite{deligne1970,deligne1972}),
on a alors un isomorphisme entre la cohomologie de De Rham 
de~$\mathscr X^*_K$ \`a p\^oles logarithmiques  le long de~$\mathscr Y_K$
et la cohomologie de De Rham de~$\mathscr X_K$:
\begin{equation}
H^i_\dR((\mathscr X^*,\mathscr Y))\otimes_W K \simeq 
H^i_\dR(\mathscr X\otimes_W K),
\end{equation}
d'o\`u un moyen concret, ne faisant intervenir
qu'un complexe de $K$-espaces vectoriels de dimensions finies,
pour calculer la cohomologie rigide de~$X$.

Sous ces hypoth\`eses,
la cohomologie log-cristalline du couple~$(X^*,Y)$
et la cohomologie de De Rham \`a p\^oles logarithmiques
$H^i_\dR((\mathscr X^*,\mathscr Y))$ sont des $W$-modules canoniquement
isomorphes. Ils d\'efinissent un $W$-r\'eseau de la cohomologie
rigide~$H^i_\rig(X/K)$, \cite{shiho2002},
munissant ainsi la cohomologie rigide 
d'une structure enti\`ere naturelle, stable par le frobenius;
cela peut \^etre utile pour les calculs effectifs. 
Dans la suite, nous qualifierons 
cette situation g\'eom\'etrique
de {\og bien relev\'ee\fg}.

La cohomologie cristalline des vari\'et\'es propres
et lisses sur~$\F$ est une cohomologie de Weil. 
D'apr\`es~\cite{katz-m74},
les espaces de cohomologie cristalline d'une telle
vari\'et\'e~$X$ ont m\^eme dimension que les espaces  correspondants
en cohomologie $\ell$-adique (pour $\ell\neq p$)
et le polyn\^ome caract\'eristique de~$F_X$ agissant
sur $H^i_\cris(X/W)\otimes K$ est le m\^eme
(avec les notations du paragraphe~\ref{subsec.l-adique})
que celui de~$\phi$ agissant sur $H^i_c(\bar X,\Q_\ell)$.
En particulier,
$\Tr(F_X|H^i_\rig(X/K))$ est un entier et il v\'erifie la majoration
\begin{equation} \label{eq.majoration.weilII}
\abs{\Tr(F_X|H^i_\rig(X/K)) } \leq q^{i/2} \dim_K H^i_\rig(X/K).
\end{equation}
D'apr\`es~\cite{chiarellotto1998}, cette majoration vaut encore
(de m\^eme bien s\^ur que sa grande cousine $\ell$-adique,
\'etablie par~\textsc{Deligne}~\cite{deligne1980})
lorsque $X$ est seulement  suppos\'e lisse sur~$\F$.

\paragraph{Principe de l'algorithme de \textsc{Kedlaya}}
Soit $X$ une vari\'et\'e alg\'ebrique d\'efinie sur~$\F$.
Pour calculer le cardinal de~$X(\F)$, \textsc{Kedlaya}
sugg\`ere de calculer les polyn\^omes caract\'eristiques du frobenius
agissant sur la cohomologie rigide de~$X$ et d'en d\'eduire $\Card{X(\F)}$
par la formule des traces de Lefschetz.

Supposons que l'on soit dans une situation g\'eom\'etrique bien relev\'ee
et que $\mathscr X$ soit un sous-sch\'ema ferm\'e
de l'espace affine de dimension~$n$,
donn\'e par des g\'en\'erateurs $(f_1,\dots,f_m)$ de son id\'eal
dans~$W[x_1,\dots,x_n]$. Notons $A$ l'anneau de~$\mathscr X$
et $A^\dagger$ son compl\'et\'e faible.
L'algorithme de \textsc{Kedlaya} est le suivant.
\begin{enumerate}
\item Calculer le rel\`evement~$\phi$ du frobenius sur~$A^\dagger$,
c'est-\`a-dire donner un algorithme pour calculer l'image
d'un \'el\'ement donn\'e \`a une pr\'ecision $p$-adique arbitraire.
Compte tenu de la condition de surconvergence impos\'ee aux
s\'eries, un tel algorithme ne manipule que des polyn\^omes;
\item Calculer des formes diff\'erentielles sur~$\mathscr X$ 
dont les classes forment une base de la cohomologie de De Rham
de~$\mathscr X_K$ (il est judicieux, mais pas n\'ecessaire,
de choisir des formes diff\'erentielles enti\`eres, 
au sens o\`u leur classe de cohomologie 
appartient \`a la cohomologie log-cristalline) ;
\item La cohomologie de Monsky--Washnitzer de~$X$ est
celle du complexe $\Omega^{*}_{A^\dagger/W}\otimes K$ 
des formes diff\'erentielles surconvergentes, mais les
classes des formes diff\'erentielles alg\'ebriques pr\'ec\'edemment calcul\'ees 
en forment une base;
Donner un algorithme {\og de r\'eduction\fg}
calculant \`a une pr\'ecision $p$-adique arbitraire
la classe de cohomologie d'une forme ferm\'ee surconvergente
(connue \`a une pr\'ecision suffisante);
\item En d\'eduire une approximation $p$-adique de la matrice~$M$
de l'endomorphisme semi-lin\'eaire~$\phi$ de~$H^k_\MW(X/K)$, 
pour $0\leq k\leq \dim X$. 
Le point d\'elicat est que la primitive d'une forme exacte \`a coefficients
enti\`ere peut faire appara\^{\i}tre des d\'enominateurs
(exemple: $x^{p-1}\,\mathrm dx$), si bien que la structure
enti\`ere de~$\Omega^1_{A^\dagger}$ n'induit pas sur $H^k_\MW(X/K)$
sa structure enti\`ere donn\'ee par la cohomologie cristalline
(lemme~2 de~\cite{kedlaya2001}, voir aussi 
le th.~2.2.5 de~\cite{abbott-kedlaya-roe2006} pour un \'enonc\'e
g\'en\'eral);
\item En d\'eduire  une approximation $p$-adique 
de la matrice
$M\sigma(M)\cdots\sigma^{d-1}(M)$,
de l'endomorphisme $K$-lin\'eaire $F_X=\phi^q$ de~$H^k_\MW(X/K)$,
puis de sa trace;
\item 
Compte tenu de l'in\'egalit\'e~\eqref{eq.majoration.weilII},
on peut en d\'eduire 
la trace elle-m\^eme, puis \'eventuellement $\Card{X(\F)}$,
si la pr\'ecision atteinte \`a l'\'etape pr\'ec\'edente est suffisante.
\end{enumerate}

\paragraph{Le cas des courbes hyperelliptiques}
La m\'ethode que nous venons d'esquisser est susceptible de s'appliquer
dans des situations tr\`es g\'en\'erales, mais une pr\'esentation d\'etaill\'ee
ne semble disponible dans la litt\'erature que pour les courbes
et le compl\'ementaire d'une hypersurface de~$\P^3$.
Nous consid\'erons ci-dessous le cas des courbes hyperelliptiques
qui faisait l'objet de l'article 
initial de~\textsc{Kedlaya}~\cite{kedlaya2001}
lorsque $p>2$ et que \textsc{Denef} 
et~\textsc{Vercauteren}~\cite{denef-vercauteren2006b}
ont \'etendu  au cas $p=2$;
voir aussi \cite{vercauteren2002,edixhoven2003}.

Soit $X^*$ une courbe hyperelliptique d\'efinie sur~$\F$.
C'est un rev\^etement double de la droite projective; les points
de ramification de ce rev\^etement sont appel\'es \emph{points de Weierstrass};
notons~$w$ leur nombre.
Si $p\neq 2$, on a $w=2g+2$, mais si $p=2$,
toute valeur de~$w$ telle que $1\leq w\leq g+1$ est possible.

Soit $Y$ l'ensemble des points de Weierstrass et posons
$X=X^*\setminus Y$.
De la cohomologie rigide de~$X$, $X^*$, $Y$, on sait un certain nombre
de choses 
\begin{itemize}
\item on a $H^0_\rig(X/K)=H^0_\rig(X^*/K)=K$,
les endomorphismes de Frobenius \'etant donn\'es par~$\sigma$ ;
\item l'espace $H^0_\rig(Y/K)$ est la somme des espaces
$H^0_\rig(P/K)$, o\`u $P$ parcourt les points de~$Y$;
en outre, si $P$ est de degr\'e~$d$ sur~$\F$, $H^0_\rig(P/K)$
s'identifie \`a l'extension non ramifi\'ee de degr\'e~$d$ de~$K$
munie de son endomorphisme de Frobenius;
\item on a $H^2_\rig(X^*/K)=K$ et le frobenius est donn\'e par~$p\sigma$.
\end{itemize}
De la suite exacte  de localisation pour la cohomologie
rigide \`a supports compacts et de la dualit\'e
de Poincar\'e,  on d\'eduit alors une suite
exacte
\begin{equation}
0 \ra H^1_\rig(X^*/K)\ra H^1_\rig(X/K) \ra H^0_\rig(Y/K)(-1)\ra H^2_\rig(X^*/K)\ra0 .
\end{equation}
En outre, les morphismes de cette suite exacte commutent aux frobenius,
le~$(-1)$ au milieu de cette formule signifiant que 
le frobenius de~$H^0_\rig(Y/K)$ est multipli\'e par~$p$.
Par ailleurs, l'involution hyperelliptique~$\eps$
agit sur ces espaces de cohomologie
et les d\'ecoupe en une partie paire et une partie impaire,
not\'ees respectivement d'un symbole~$+$ et~$-$, donnant lieu \`a une
suite exacte
\begin{equation}\label{eq.coh.affine}
0\ra  H^1_\rig(X/K)^+ \ra H^0_\rig(Y/K)(-1)
   \ra H^2_\rig(X^*/K) \ra 0
\end{equation}
(car $H^1_\rig(X^*/K)^+ \simeq H^1_\rig(\P^1/K)=0$)
et un isomorphisme
\begin{equation}
 H^1_\rig(X^*/K)^- \simeq H^1_\rig(X/K) .
\end{equation}
Il s'agit donc de calculer $H^1_\rig(X/K)^-$.

Les r\'ef\'erences ci-dessus ne traitent en fait que le cas
o\`u l'un des points de Weierstrass est rationnel.
Dans ce cas, la courbe~$X^*$ 
poss\`ede une \'equation  plane  de la forme (affine)
\begin{subequations}
\begin{equation}\label{eq.hyperell.impair}
y^2=f(x), \qquad f\in\F[x], \quad \deg(f)=2g+1
\end{equation}
si $p\neq 2$; 
lorsque $p=2$, elle a une \'equation du type
\begin{equation}\label{eq.hyperell.p=2}
y^2+f(x)y=g(x), \qquad \deg(f)\leq g, \quad \deg(g)=2g+1.
\end{equation}
\end{subequations}
Avec ces \'equations, l'involution hyperelliptique
est donn\'ee par \mbox{$(x,y)\mapsto (x,-y)$}, resp. $(x,y)\mapsto (x,-f(x)-y)$;
la projection vers~$\P^1$ est donn\'ee par l'application \mbox{$(x,y)\mapsto x$}; 
elle applique le point de Weierstrass choisi sur le point
\`a l'infini de~$\P^1$.  
Les autres points de Weierstrass sont ceux de coordonn\'ees~$(x,y)$
tels que $f(x)=0$.
Lorsque $p\neq 2$, ce sont ceux d'ordonn\'ee nulle;
lorsque $p=2$, 
l'application d'une transformation convenable
de la forme $(x,y)\mapsto (x,y+\alpha(x))$ 
permet de supposer que $g$ s'annule en chacun de ces points;
leur ordonn\'ee est encore nulle. Soit $h$ le produit
des facteurs irr\'eductibles de~$f$; le diviseur de~$h$
est \'etale sur~$\F$ et, au point \`a l'infini pr\`es,
a pour support les points de Weierstrass de~$X^*$.

Soit $\bar B$ l'anneau $\F[x,h(x)^{-1}]$ ; l'anneau~$\bar A$ 
de la courbe affine~$X$
est un $\bar B$-module libre de rang~$2$, de base~$(1,y)$, $y$
v\'erifiant l'\'equation~\eqref{eq.hyperell.impair} si $p\neq 2$, 
resp.~\eqref{eq.hyperell.p=2} si $p=2$.
En particulier, $\bar A$ est \'etale sur~$\bar B$.

Lorsque $p\neq 2$, choisissons un rel\`evement~$\tilde f$
de~$f$ dans l'anneau~$W[x]$ de m\^eme degr\'e que~$f$
et posons $\tilde h=\tilde f$.
L'\'equation analogue \`a~\eqref{eq.hyperell.impair} d\'efinit une courbe
hyperelliptique~$\mathscr X^*$ de genre~$g$ 
qui rel\`eve la courbe~$X^*$;
le sch\'ema~$\mathscr Y$ des points de Weierstrass de~$\mathscr X^*$
est \'etale sur~$W$ et la $W$-courbe affine
$\mathscr X=\mathscr X^*\setminus\mathscr Y$ rel\`eve~$X$.
Notons~$B$ l'anneau~$W[x,\tilde f(x)^{-1}]$;
l'anneau~$A$ de~$\mathscr X$ est \'egal \`a \mbox{$B[y]/(y^2-\tilde f)$}
et est \'etale sur~$B$.

Lorsque $p=2$, le diviseur des points de Weierstrass d'une
courbe~$\mathscr X^*$ qui rel\`eve~$X^*$ n'est jamais \'etale sur~$W$;
il convient alors de choisir un point de Weierstrass par classe r\'esiduelle. 
De mani\`ere pr\'ecise,
\textsc{Denef} et~\textsc{Vercauteren} commencent par
relever~$h$ en un polyn\^ome~$\tilde h$ de m\^eme degr\'e
puis exigent que chaque facteur irr\'eductible de~$\tilde h$
divise~$\tilde f$ et~$\tilde g$ avec la m\^eme multiplicit\'e
que celle dont le facteur irr\'eductible correspondant 
de~$h$ divise respectivement~$f$ et~$g$.
En particulier $\tilde g$ s'annule en toute racine de~$\tilde f$.
L'\'equation analogue \`a~\eqref{eq.hyperell.p=2}
d\'efinit alors une courbe hyperelliptique~$\mathscr X^*$ de genre~$g$
qui rel\`eve la courbe~$X^*$; notons~$\mathscr Y$ 
la r\'eunion du point \`a l'infini 
et du lieu d\'efini par le polyn\^ome~$\tilde h$.
C'est un sous-sch\'ema \'etale de~$\mathscr X^*$ dont la r\'eduction
modulo~$p$ a pour support 
l'ensemble des points de Weierstrass de~$\mathscr X^*$.
Posons $\mathscr X=\mathscr X^*\setminus\mathscr Y$.
Posons $B=W[x,\tilde h(x)^{-1}]$; l'anneau~$A$ de~$\mathscr X$
est $B[y]/(y^2+\tilde f y-\tilde g)$.

Dans les deux cas,
l'anneau~$B^\dagger$ peut \^etre d\'ecrit explicitement,
comme un anneau de s\'eries en~$x$ et~$\tilde h(x)^{-1}$
dont les coefficients tendent assez vite vers~$0$.
L'anneau~$A^\dagger$ est alors une $B^\dagger$-alg\`ebre
\'etale, libre de rang~$2$, de base~$(1,y)$.
Pour $p=2$, ce sont les m\^emes anneaux que ceux qu'on aurait obtenu
en rempla\c{c}ant~$\mathscr Y$ par le sous-sch\'ema des points de Weierstrass.

Pour relever le frobenius, nous commen\c{c}ons par choisir
sur l'anneau~$W[x]$ l'unique rel\`evement $\sigma$-lin\'eaire
tel que $x\mapsto x^p$, autrement dit $\sum a_n x^n\mapsto \sum \sigma(a_n)x^{pn}$.
Comme $A$ est \'etale sur~$W[x]$, cet homomorphisme
s'\'etend de mani\`ere unique en un homomorphisme $\sigma$-lin\'eaire, not\'e~$\phi$,
de~$A^\dagger$ dans lui-m\^eme.
M\^eme si cela r\'esulte 
d'une variante due \`a~\textsc{Bosch}~\cite{bosch1981},
du th\'eor\`eme d'approximation
d'\textsc{Artin} pour les anneaux de s\'eries
surconvergentes,
nous devons calculer~$\phi$ explicitement,
et en particulier d\'eterminer l'\'el\'ement~$\phi(y)$ de~$A^\dagger$
qui rel\`eve~$y^p$ et tel que
$\phi(y)^2=\tilde f^\sigma (x^p)$ lorsque $p\neq 2$
et $\phi(y)^2+\tilde f^\sigma(x^p) \phi(y)-\tilde g^\sigma(x^p)$
si $p=2$.
Supposons d'abord $p\neq 2$;
il existe  
un \'el\'ement $h\in W[x]$ tel que 
$\tilde f^\sigma(x^p)=(\tilde f(x))^p+ph(x)$, 
car $\sigma$ rel\`eve l'automorphisme de Frobenius de~$\F$.
On pose alors
\begin{equation}
\phi(y)=y^p \left( 1+ p \frac{h(x)}{\tilde f(x)^p} \right)^{1/2}
= y^p \sum_{k=0}^\infty \binom{1/2}k p^k \frac{h(x)^k}{\tilde f(x)^{pk}}.
\end{equation}
La pr\'esence des coefficients~$p^k$ permet d'analyser tr\`es
facilement la convergence de cette s\'erie. 
En particulier, l'\'el\'ement $\phi(y)$
\'ecrit appartient \`a~$A^\dagger$ et v\'erifie l'\'equation consid\'er\'ee.
Dans le cas $p=2$, le principe est similaire :
comme on a pris soin de placer les points de Weierstrass sur
l'axe~$y=0$, $y$ est inversible dans~$A^\dagger$ 
et on peut chercher $\phi(y)$ sous la forme $y^pu$, avec $u\equiv 1\pmod p$ ;
je renvoie \`a~\cite{denef-vercauteren2006b}, lemma~1,
pour les d\'etails.
Signalons aussi qu'en pratique,
$\phi(y)$ n'est pas calcul\'e en d\'eveloppant des
s\'eries enti\`eres mais en appliquant la m\'ethode de Newton.

Comme $\tilde Y$ est \'etale sur~$W$,
le th\'eor\`eme de comparaison
entre cohomologies de De Rham alg\'ebrique et cohomologie rigide
\'evoqu\'e plus haut entra\^{\i}ne que
$H^1_\MW(X/K)$  s'identifie au premier groupe de cohomologie
de De Rham de~$X\otimes K$.\footnote{%
C'est l\`a qu'intervient le choix de~$\mathscr Y$ en caract\'eristique~$2$:
il ne fallait pas enlever tous les points de Weierstrass de la courbe
hyperelliptique~$\mathscr X$, mais seulement un par point
de Weierstrass de la fibre sp\'eciale.}
Il en est de m\^eme de la partie~$-$,
ce qui montre que $H^1_\MW(X/K)^-$ admet
pour base les classes $[x^i y^{-1}\mathrm dx]$, pour $0\leq i\leq 2g-1$.
Comme l'algorithme final n'utilisera qu'une pr\'ecision finie,
il est de toutes fa\c{c}ons n\'ecessaire 
de contr\^oler la valeur absolue des coefficients des diff\'erentielles
exactes mises en jeu par une telle identification.
L'{\og algorithme de r\'eduction\fg}
permet d'\'ecrire une forme diff\'erentielle ferm\'ee surconvergente
\`a coefficients entiers~$\omega$ comme somme de deux termes:
d'une part une combinaison lin\'eaire explicite
des formes diff\'erentielles qui constituent 
la base de~$H^1_\MW(X/K)$,
d'autre part une forme reste dont un multiple explicite
est la diff\'erentielle d'une forme diff\'erentielle surconvergente 
\`a coefficients entiers.

Supposons ces deux points acquis et soit $i\in\{0,\dots,2g-1\}$.
Pour calculer l'image par~$\phi$ de la classe
de la forme diff\'erentielle~$\omega_i=x^iy^{-1}\mathrm dx$,
il reste \`a effectuer les calculs suivants:
\begin{itemize}
\item d\'evelopper en s\'erie la forme diff\'erentielle
\begin{equation}
\begin{split}
\phi(\omega_i) &= \phi(x)^i\phi(y)^{-1}\mathrm d\phi(x)
= px^{pi+p-1}\phi(y)^{-1}\mathrm dx \\
& = px^{pi+p-1}\tilde f(x)^{-(p-1)/2)}\left( 1+ p \frac{h(x)}{\tilde f(x)^p} \right)^{1/2} y^{-1}\mathrm dx,\end{split}
\end{equation}
\`a une pr\'ecision $p$-adique suffisante ;
\item \'ecrire le terme principal comme une somme
de deux termes: le premier est une combinaison
lin\'eaire des formes diff\'erentielles $\omega_k$,
pour $0\leq k\leq 2g-1$,  le second (qu'en fait
on n'\'ecrit pas) est une forme diff\'erentielle 
dont un multiple explicite est une forme exacte;
\item de m\^eme, l'image du terme reste dans la cohomologie
sera alors  combinaison lin\'eaire
des classes~$\omega_i$ avec des coefficients $p$-adiquement petits.
\end{itemize}
Si la pr\'ecision a \'et\'e choisie assez grande, on en d\'eduit
une approximation
de l'image de~$[\omega_i]$ par~$\phi$, 
donc une approximation de la matrice~$M$ de~$\phi$.
La matrice de~$F_X=\phi^d$ est, quant \`a elle,
donn\'ee par $M\sigma(M)\dots\sigma^{d-1}(M)$.
Pour finir, si la pr\'ecision est suffisante, on peut calculer 
la trace de~$F_X$ (donc le cardinal de~$X(\F)$) et
le d\'eterminant de~$1-tF_X$ (donc la fonction z\^eta de~$X$).

\paragraph{Complexit\'e et g\'en\'eralisations}
La complexit\'e de l'algorithme que nous avons grossi\`erement
d\'ecrit d\'epend de la pr\'ecision requise pour effectuer
les calculs; je renvoie aux articles cit\'es ainsi
qu'\`a~\cite{gaudry-gurel2003} pour l'analyse de cette complexit\'e.
Lorsque $p\neq 2$, il en ressort qu'elle est~$\Otilde (pg^4d^3)$ en temps
et~$\Otilde(pg^3d^3)$ en espace.  (Rappelons que $q=p^d$.)
Lorsque $p=2$, la complexit\'e est un peu moins bonne en temps,
\`a savoir $\Otilde(g^5d^3)$. Tout r\'ecemment, des
id\'ees remontant aux \textsc{Chudnovsky} ont permis
de faire baisser la d\'ependance en~$p$ de lin\'eaire \`a~$\sqrt p$
(voir~\cite{harvey2006}, ainsi que \cite{bostan-gaudry-schost2007}).

En pratique, l'algorithme de~\textsc{Kedlaya}
a permis de calculer le cardinal de courbes
hyperelliptiques de genres~$\leq 4$ en quelques
minutes; le produit~$gd\log_2 p$ (approximativement
\'egal au logarithme en base~$2$ du cardinal de la jacobienne) 
\'etant de l'ordre de~$200$
($p=2$, \cite{denef-vercauteren2006b};
$p=251$, \cite{gaudry-gurel2003}).

Par ailleurs, il a donn\'e lieu \`a un certain nombre de g\'en\'eralisations:
courbes superelliptiques ($y^m=f(x)$, \cite{gaudry-gurel2001}),
courbes~$C_{a,b}$ (rev\^etements de~$\P^1$ totalement ramifi\'es
\`a l'infini) dans~\cite{denef-vercauteren2006}, etc.
Pour les courbes, l'algorithme le plus g\'en\'eral semble
celui de \textsc{Castryck}, \textsc{Denef} et \textsc{Vercauteren}
dans~\cite{castryck-denef-vercauteren2006} qui concerne
les courbes planes
qui sont non d\'eg\'en\'er\'ees par rapport \`a leur polytope de Newton;
si $p$ est fix\'e,
la complexit\'e en temps de cet algorithme est~$\Otilde(d^3g^{6{,}5})$,
celle en espace est~$\Otilde(d^3g^4)$.

Enfin, l'article~\cite{abbott-kedlaya-roe2006}
utilise cette m\'ethode de calcul de la cohomologie $p$-adique
pour \'evaluer le rang du nombre de Picard
de surfaces projectives lisses,
ces deux quantit\'es  \'etant reli\'ees par
une conjecture de~\textsc{Tate}
qui relie le rang du groupe de N\'eron-Severi d'une surface
projective lisse~$S$ d\'efinie sur~$\F$
\`a la multiplicit\'e de la valeur propre~$q$ de l'endomorphisme~$F_S$
sur $H^2_\rig(X/K)$.
(Pour d\'emontrer une majoration du nombre de Picard,
la conjecture de~\textsc{Tate} n'est bien s\^ur pas n\'ecessaire.)

\subsection{Variation de la cohomologie $p$-adique}

La complexit\'e en temps des m\'ethodes $p$-adiques d\'ecrites jusqu'ici 
est toujours exponentielle en la dimension de l'espace ambiant.
Le dernier algorithme de ce texte, introduit par~\textsc{Lauder}~\cite{lauder2004} en~2002, vise \`a supprimer ce d\'efaut
en tirant parti de la variation de la cohomologie $p$-adique
dans une famille, longuement \'etudi\'ee par~\textsc{Dwork}
dans les ann\'ees~60.

Plusieurs incarnations de cet algorithme
sont actuellement disponibles:
l'article original de~\textsc{Lauder}~\cite{lauder2004},
r\'edig\'e dans le cadre de la th\'eorie de~\textsc{Dwork},
concerne les familles d'hypersurfaces  dont un membre
est une hypersurface diagonale.
Ind\'ependamment, N.~\textsc{Tsuzuki}~\cite{tsuzuki2003b}
avait propos\'e un algorithme qui calcule
des sommes de Kloosterman et le nombre de points des
rev\^etements d'Artin-Scheier de~$\gm$.
R.~\textsc{Gerkmann}~\cite{gerkmann2005}
dans le cas des courbes elliptiques,
puis H.~\textsc{Hubrechts} 
(voir \cite{hubrechts2006a,hubrechts2006b, hubrechts2006c}
ainsi que l'esquisse~\cite{lauder2005})
pour les courbes hyperelliptiques,
ont utilis\'e cette m\'ethode
et propos\'e un algorithme calculant le cardinal d'une telle courbe,
de genre~$g$, d\'efinie sur un corps fini de cardinal~$p^d$,
dont la complexit\'e en espace est~$\Otilde(pd^2g^4)$;
cela am\'eliore l'algorithme de~\textsc{Kedlaya} par rapport
au param\`etre~$d$.
Plus r\'ecemment, \textsc{Lauder}~\cite{lauder2006} a appliqu\'e 
cette m\'ethode \`a un pinceau de Lefschetz 
dont la vari\'et\'e est l'espace total, obtenant ainsi un algorithme
r\'ecursif pour calculer la fonction z\^eta d'une vari\'et\'e
alg\'ebrique sur un corps fini.
Ces derniers algorithmes sont r\'edig\'es en termes de cohomologie rigide
ou Monsky--Washnitzer.

Comme cette \'enum\'eration le montre,
la port\'ee de cette m\'ethode est tr\`es g\'en\'erale;
je me contente cependant d'en esquisser ici le principe 
dans le cas d'une famille de courbes.
On reprend les notations $\F$, $q$, $p$, $W$ et~$K$.

\paragraph{La d\'eformation}
Le but de l'algorithme est de calculer la fonction z\^eta d'une courbe~$X_1$,
ou plus simplement le cardinal de~$X_t(\F)$, lorsque $X_t$
est la fibre en un point~$t$ d'une 
famille de courbes param\'etr\'ee par un ouvert~$U$ de~$\A^1$.
On peut imaginer par exemple une famille de courbes
hyperelliptiques donn\'ee par une \'equation hyperelliptique
\[ X_t\colon  y^2 = f_t(x), \]
o\`u $f_t\in \F[x,t]$ est de degr\'e~$2g+1$ en~$x$, $\F$ \'etant un corps fini
de cardinal~$q$. Dans ce cas, on notera $U$ le plus grand ouvert de~$\A^1$
au-dessus duquel $f_t$ est s\'eparable et de degr\'e~$2g+1$; son compl\'ementaire
est le lieu d\'efini par le discriminant~$D(t)$ de~$f_t$.
Si l'on s'int\'eresse \`a des courbes affines, il est prudent de
consid\'erer une famille propre et lisse $u\colon X^*\ra U$,
un diviseur $Y$ \'etale sur~$U$, puis de poser $X=X^*\setminus Y$.
Pour $t\in U$, $X^*_t$, $Y_t$, $X_t$  d\'esignent les fibres de~$X^*$, $Y$,
$X$ au-dessus de~$t$.

L'algorithme de d\'eformation repose sur le fait que lorsque~$t$
varie (en un certain sens), les matrices de Frobenius de~$H^1_\rig(X_t/K)$
v\'erifient une \'equation diff\'erentielle. L'id\'ee de~\textsc{Lauder}
est de partir de la matrice en un point $t=0$, suppos\'ee
d\'ej\`a calcul\'ee par une autre m\'ethode, et de r\'esoudre 
cette \'equation diff\'erentielle pour arriver au point~$t$.

Il y a en cohomologie rigide une notion d'image directe
et les $K$-espaces vectoriels de cohomologie rigide $H^1_\rig(X_t,K)$,
pour $t\in U$, sont (en un certain sens que nous pr\'eciserons
plus bas) les fibres d'un objet $R^1(u_\rig)_*(X/\P^1)$
(cohomologie rigide relative)
qui est, au moins conjecturalement\footnote{
Avec les hypoth\`eses impr\'ecises ci-dessus, c'est vraisemblable,
mais je ne l'ai pas v\'erifi\'e en d\'etail;
dans le cas, pr\'esent\'e plus bas, des courbes hyperelliptiques, cela r\'esulte
de l'\'etude explicite faite par~\textsc{Hubrechts}.},
un $F$-isocristal surconvergent sur~$U$.
Pour en donner une description concr\`ete, on suppose
toute la situation bien relev\'ee en caract\'eristique~$0$
en se donnant un ouvert~$\mathscr U$ de~$\P^1_W$, compl\'ementaire d'un 
$W$ sch\'ema \'etale d'\'equation $\tilde D(t)=0$,
une famille propre et lisse $u\colon \mathscr X^*\ra\mathscr U$
et un diviseur $\mathscr Y$ de~$\mathscr X^*$, \'etale sur~$\mathscr U$.

Notons $S$ l'anneau $W[t,\tilde D(t)^{-1}]$ et $S^\dagger$ son compl\'et\'e faible;
soit $\sigma$ l'endomorphisme semi-lin\'eaire de~$S$ qui applique $t$ sur~$t^p$.
Identifions les \'el\'ements de~$S^\dagger_K$ \`a des s\'eries
formelles en~$t$; alors, si  $f\in S^\dagger_K$,  $\sigma(f(t))=f^\sigma (t^p)$,
o\`u $f^\sigma$ est la s\'erie formelle obtenue en appliquant~$\sigma$
aux coefficients de~$f$.  

Alors, $R^1(u_\rig)_*(X/\P^1)$ est un $S^\dagger_K$-module $E$,
libre de rang~$m=2g-1+r$ (o\`u $g$ est le genre de~$X_t^*$ et $r$ le
nombre de points g\'eom\'etriques de~$Y_t$, c'est la version
en famille de~\eqref{eq.coh.affine}),
muni d'une connexion 
\begin{equation}
  \nabla\colon
        E\ra  E\otimes_{S^\dagger} \Omega ^1_{S^\dagger_K}
\end{equation}
et d'un isomorphisme horizontal $\phi\colon \sigma^*E\simeq E$.
En outre, pour tout $t\in \F$, de rel\`evement de Teichm\"uller $\tilde t\in W$
consid\'er\'e comme un homomorphisme d'anneaux $S^\dagger_K\ra K$,
$E\otimes_{\tilde t} K$ s'identifie \`a $H^1_\rig(X_t/K)$
et $\phi\otimes_{\tilde t} 1$ s'identifie \`a l'endomorphisme de Frobenius
de la cohomologie rigide (noter que $\sigma(\tilde t)=\tilde t^p$).
La connexion~$\nabla$ est bien s\^ur
un avatar de la connexion de Gauss-Manin en cohomologie de De Rham
(et devrait \^etre la m\^eme lorsqu'on dispose d'un th\'eor\`eme 
de comparaison d'images directes entre cohomologie de De Rham
et cohomologie rigide).

Fixons une base $(e_1,\dots,e_m)$ de~$E$ et soit $G(t) \dt$ 
la matrice de~$\nabla$
dans cette base; on a donc, pour tout vecteur colonne~$X(t)$ 
\`a coefficients dans~$S^\dagger_K$,
\begin{equation}
\nabla \big((e_1,\dots,e_m)X(t) \big) = (e_1,\dots,e_m)
(G(t)X(t)+\frac{\mathrm dX(t)}{\dt}) \, \mathrm dt.
\end{equation}
Soit aussi $F(t)$ la matrice de~$\phi$ dans les
bases $(e_i\otimes 1)$ et~$(e_i)$; on a ainsi
\begin{equation}
\phi \big( (e_1,\dots,e_m) X(t)\otimes 1 \big)  = (e_1,\dots,e_m)
           (F(t) X^\sigma (t^p))
\end{equation}
et l'horizontalit\'e de~$\phi$ se traduit par l'\'equation diff\'erentielle
\begin{equation}\label{eqdiff.F-nabla}
F'(t) + G(t) F(t) = F(t) G^\sigma (t^p) p t^{p-1}.
\end{equation}

\paragraph{R\'esolution et surconvergence}
Soit $X(t)$ une matrice fondamentale de l'\'equation diff\'erentielle homog\`ene
\begin{equation}\label{eqdiff.homogene}
X'(t)+G(t)X(t)=0,  \qquad X(0)=\mathrm I_m.
\end{equation}
Ses coefficients sont des s\'eries formelles \`a coefficients dans~$K$.
La d\'emonstration du th\'eor\`eme de Cauchy fait appara\^{\i}tre
des d\'enominateurs ; par exemple, le rayon de convergence $p$-adique de
la s\'erie exponentielle, solution de $y'=y$
n'est que~$p^{-1/(p-1)}$. Toutefois, une astuce due \`a~\textsc{Dwork}
reposant sur la formule~\eqref{eq.frobenius.dwork} ci-dessous montre
que le rayon de convergence de $X(t)$ est \'egal \`a~$1$, de sorte que
$X(t)$ converge dans le disque unit\'e \emph{ouvert}.

Par horizontalit\'e de~$\phi$, l'image 
de~$(e_1,\dots,e_m)X(t)\otimes 1$ par~$\phi$ est 
annul\'ee par~$\nabla$; il existe donc une matrice \emph{constante}~$C$
telle que l'on ait
$      F(t) X^\sigma(t^p) = X(t) C$,
d'o\`u $F(0)=C$ et
\begin{equation}\label{eq.frobenius.dwork}
 F(t) = X(t) F(0) X^\sigma(t^p)^{-1}.
\end{equation}
Ayant calcul\'e~$X(t)$, on peut tout autant calculer~$F(t)$,
pourvu que la matrice~$F(0)$ ait \'et\'e calcul\'ee auparavant,
c'est-\`a-dire que l'on connaisse l'action du frobenius
sur la cohomologie du membre~$X_0$ de notre famille.
Ce n'est  a priori qu'une s\'erie formelle \`a coefficients dans~$K$,
de m\^eme rayon de convergence que~$X(t)$.

Mais $F(t)$ est bien plus : c'est une matrice
\`a coefficients dans~$S^\dagger_K$. Ce fait g\'en\'eral,
qui traduit la {\og surconvergence\fg} de l'isocristal
$R^1(u_\rig)_*(X/\P^1)$ \'etait d\'ej\`a apparu dans cet expos\'e:
la s\'erie~$\theta$ de l'\'equation~\eqref{eq.theta.dwork}
dont on avait mentionn\'e la d\'ecroissance vers~$0$ des coefficients
est celle qui correspond \`a l'isocristal de Dwork sur la droite affine
(voir par exemple~\cite{berthelot1988}, p.~28-29).

Puisque le d\'eveloppement en s\'erie formelle en l'origine
d\'efinit un homomorphisme \emph{injectif} de~$S^\dagger_K$ dans~$K[[t]]$,
on peut exprimer de mani\`ere unique chaque entr\'ee de
la matrice~$F(t)$ comme
une s\'erie de la forme
\begin{equation}\label{eq.fn}
   f(t) = \sum _{n=-\infty}^\infty f_n \tilde D(t)^n,
\end{equation}
o\`u les $f_n$ sont des polyn\^omes de degr\'es~$<\deg\tilde D$
dont les coefficients tendent rapidement vers~$0$ lorsque $\abs n\ra\infty$.

Supposons pour l'instant que l'on ait r\'eussi \`a \'ecrire les entr\'ees
de la matrice~$F(t)$ sous cette forme.
Si $t\in\F$ est un point de~$U$,
on peut alors \'evaluer $F$ en le rel\`evement de Teichm\"uller~$\tilde t$
de~$t$, car $\abs{t}=\abs{\tilde D(t)}=1$, d'o\`u une matrice $F(\tilde t)$
qui est celle du frobenius agissant sur la cohomologie rigide de
la fibre~$X_t$.

Il reste \`a expliquer comment l'on peut calculer les~$f_n$.
Dans le cas sp\'ecifique d'une famille de courbes hyperelliptiques,
\textsc{Hubrechts} fait dans~\cite{hubrechts2006b}
toute l'analyse pr\'ec\'edente en grand d\'etail, et de mani\`ere explicite;
il prouve de fait la surconvergence de l'isocristal $R^1(u_\rig)_*(X/\P^1)$
en \'etablissant une minoration explicite des 
valuations $p$-adiques des coefficients des polyn\^omes~$f_n$.
Il peut alors consid\'erer l'\'equation~\eqref{eq.fn} comme
un syst\`eme d'\'equations lin\'eaires \`a coefficients dans~$K$.
En chassant les d\'enominateurs et en se limitant \`a une pr\'ecision $p$-adique
donn\'ee,
on obtient un syst\`eme lin\'eaire en dimension finie dont on peut
calculer une solution; la pr\'ecision obtenue est inf\'erieure \`a la pr\'ec\'edente
car le syst\`eme n'est pas inversible modulo~$p$.

Ainsi, \textsc{Hubrechts} est en mesure de pr\'evoir
quelle pr\'ecision initiale est n\'ecessaire dans tout ce calcul
pour, \emph{in fine}, obtenir les entr\'ees de la matrice~$F(t)$
\`a une pr\'ecision suffisante pour que la congruence qui en r\'esultera
sur la fonction z\^eta permette de la d\'eterminer.

En pratique, plut\^ot que de calculer une solution fondamentale~$X(t)$
de~\eqref{eqdiff.homogene}, puis $G(t)$ par la formule~\eqref{eq.frobenius.dwork},
\textsc{Hubrechts} donne des algorithmes it\'eratifs efficaces
pour r\'esoudre directement des \'equations diff\'erentielles
du type~\eqref{eqdiff.F-nabla}.

Pour les courbes elliptiques, l'algorithme ainsi d\'ecrit
est moins rapide que l'\textsc{agm}, ce dernier \'etant cependant restreint
en pratique au cas $p=2$.
Selon les donn\'ees pr\'esent\'ees \`a la fin de~\cite{hubrechts2006b},
l'algorithme est plus rapide que l'algorithme~\textsc{sea}
lorsque le degr\'e est au moins~$100$ si $p=3$, et $40$ si $p=7$.
Le d\'enombrement d'une courbe elliptique sur un corps
de cardinal~$3^{500}$ n\'ecessite un peu plus d'une heure de calcul
et 50~Mo d'espace disque.
En genre~$2$, il permit de calculer le nombre de points d'une
courbe sur un corps de cardinal~$3^{400}$ en une vingtaine d'heures
et 120~Mo d'espace disque; un tel calcul requerrait plusieurs~Go
avec l'algorithme de~\textsc{Kedlaya}.

\def\bibliofont{\def\baselinestretch{1.05}\small}
\bibliographystyle{smfplain}
\bibliography{aclab,fast,acl}
\addressindent.3\textwidth
\end{document}